\documentclass[12pt,a4paper]{article}
\usepackage{amsmath,amssymb,latexsym,amsthm}
\usepackage[english]{babel}

\begin{document}

\newtheorem*{Dfn}{Definition}
\newtheorem{Theo}{Theorem}
\newtheorem{Lemma}[Theo]{Lemma}
\newtheorem{Prop}[Theo]{Proposition}
\newcommand{\Pro}{\noindent{\em Proof. }}
\newcommand{\Rem}{\noindent{\em Remark. }}

\title{$3$-dimensional Bol loops as sections in non-solvable Lie groups}
\author{\'Agota Figula}
\date{}
\maketitle

\begin{abstract} Our aim in this paper is to classify the $3$-dimensional connected 
differentiable
 global Bol loops, which have a non-solvable group as the group topologically 
generated by their left translations
and to describe their relations to metric space geometries.
The classification of global differentiable Bol loops significantly differs  
from the classification of local differentiable Bol loops. 
We treat the differentiable Bol loops as images of global differentiable 
sections $\sigma :G/H \to G$  
such that for all $r,s \in \sigma (G/H)$ the element $rsr$ lies in 
$\sigma (G/H)$, where $H$ is the stabilizer of the identity $e$ of $L$ in $G$. 
\end{abstract}



\section{Introduction}

Differentiable diassociative loops can be treated following the ideas of 
Sophus Lie since for their tangential objects, binary Lie algebras, Lie's 
first theorem  
and  Lie's third theorem are satisfied: To every binary Lie algebra there  
exists (up to local 
isomorphisms) a unique differentiable  diassociative local loop. Kuzmin, Kerdman,  
and Nagy have proved that any local 
differentiable Moufang loop (these loops form a subclass of the class of  
diassociative loops) can uniquely
be  embedded into a connected simply connected global one 
(cf. \cite{kerdman}, \cite{kuzmin}, \cite{nagy}). This 
problem is more difficult for diassociative loops. Namely, there are local 
differentiable diassociative loops such that the corresponding global 
analytic loops do not exist (cf. \cite{grishkov2}). In  
 \cite{grishkov2} the author gives in special cases the exact condition  
when a local 
diassociative loop can be  embedded in a connected simply connected 
global one.  
\newline
The present research on loops turns to  
such classes $\cal C$ of differentiable loops which have local forms  
determined in a unique way by algebras the properties of which are deduced 
by the 
properties of the tangential objects of loops in the class $\cal C$.
The most important class $\cal C$ of such loops are the Bol loops. Their 
tangential objects, the Bol algebras, can be seen as  Lie triple systems 
with an additional binary operation (cf. \cite{loops} pp. 84-86, Def. 6.10.).  
As known the Lie triple systems are in one-to-one correspondence to (global) 
simply connected symmetric spaces (cf. \cite{symmetric}, \cite{loops} Sec. 6).
Hence there is a strong connection between the theory of differentiable Bol 
loops and the theory of symmetric spaces.
In particular the theory of  connected differentiable Bruck loops 
( which form a 
subclass of the class of Bol loops) is essentially the theory of affine 
symmetric spaces (cf. \cite{loops} Sec. 11).   
\newline
The $1$-dimensional connected differentiable Bol loops are isomorphic either 
to the group $\mathbb R$ or to $SO_2(\mathbb R)$. 
The smallest connected differentiable proper Bol loops are realized on   
$2$-dimensional manifolds. 
There exist precisely two isotopism classes of  proper
$2$-dimensional connected differentiable global Bol loops. As a representative of the one class we may choose 
the hyperbolic plane loop (cf. \cite{loops}, Section 22);  this loop is the unique Bruck loop in the isotopism class of loops, having the connected 
component of the group of hyperbolic motions as the group topologically generated by their left translations. 
As a representative of the other isotopism class may be chosen the 
$2$-dimensional Bruck loop which is realized on the pseudo-euclidean affine 
plane such that its left translations group 
is the connected component of the group of pseudo-euclidean 
motions and the elements of $L$ are the lines of positive slope in the 
pseudo-euclidean affine plane (cf. \cite{loops} Section 25). 
\newline
Our aim in this paper is to classify the $3$-dimensional connected 
differentiable
 global Bol loops, which have a non-solvable group as the group topologically 
generated by their left translations
and to describe their relations to metric space geometries.
The classification of global differentiable Bol loops significantly differs  
from the classification of local differentiable Bol loops. 
There exist much more local differentiable Bol loops than global 
differentiable Bol loops, as we will also show  in this paper.     
\newline
We treat the differentiable Bol loops as images of global differentiable 
sections $\sigma :G/H \to G$  
such that for all $r,s \in \sigma (G/H)$ the element $rsr$ lies in 
$\sigma (G/H)$, where $H$ is the stabilizer of the identity $e$ of $L$ in $G$. 
Using the present theory of differentiable Bol loops it is not difficult to 
prove that $G$ is four, five or sixdimensional. The images of sections 
for local Bol loops can be realized as the exponential images of 
Lie triple systems. But then we have to discuss which of these local Bol loops can be extended to global Bol loops.  
\newline
The results of our paper can summarized in the following 
\newline
\newline
{\bf Theorem} \begin{em} 
 There are precisely two isotopism 
classes ${\cal C}_i,\ (i=1,2)$ of the $3$-dimensional
connected differentiable simple proper Bol loops $L$ such that the group $G$ 
generated by 
the left translations $\{ \lambda_x; x \in L \}$ is a non-solvable Lie group.

The  class ${\cal C}_1$ consists of Bol loops having the simple 
Lie group $G= PSL_2(\mathbb C)$ as the group topologically generated by their left translations and as the stabilizer $H$ of 
$e \in L$ in $G$ the group $SO_3(\mathbb R)$. 
\newline
Any loop in the class ${\cal C}_1$ can be characterized by a real 
parameter $a$, where $-1 < a < 1$.
The loops $L_a$  are isotopic to the hyperbolic 
space loop
$L_0$, which is realized on the hyperbolic space by the multiplication 
$x \cdot y= \tau_{e,x} (y)$, where $\tau_{e,x}$ is the hyperbolic translation 
moving $e$ onto $x$.    
The loop $L_0$  is the only  Bruck  loop in ${\cal C}_1$. 
The loops $L_a$ and $L_b$ in ${\cal C}_1$ are isomorphic if and 
only if the angles between the tangent $3$-space ${\bf m_a}$, respectively 
${\bf m_b}$ and the tangent 
$3$-space ${\bf m}_0$ are the same 
with respect to the Cartan-Killing form.  
This is the case if and only if $b= \pm a$.
\newline
The other class ${\cal C}_2$ of simple Bol loops  consists of  
$3$-dimensional connected differentiable Bol 
loops such that the group  $G$ topologically generated by their left translations is the semidirect product 
$PSL_2(\mathbb R) \ltimes \mathbb R^3$,  where 
the action of $PSL_2(\mathbb R)$ on $\mathbb R^3$ is the adjoint action of 
$PSL_2(\mathbb R)$ on its Lie algebra, and for the stabilizer $H$ of $e \in L$ 
we have 
\newline
\centerline{ $ \left \{ \left(\pm \left( \begin{array}{rr}
\cos t & \sin t \\
-\sin t & \cos t \end{array} \right) , \left( \begin{array}{rr}
-x & y \\
y & x \end{array} \right) \right); t \in [0 , 2 \pi ],x,y \in \mathbb R 
\right \}$.}  
\newline
Any loop in this class can be characterized by three real parameters 
$a,b,c$, with $a^2+b^2<1$. 
In ${\cal C}_2$ there is only one isomorphism class of Bruck loops, which 
consists 
of the loops $\hat{L}_c=L_{0,0,c}$. As a representative of this isomorphism 
class may be chosen the pseudo-euclidean space loop 
$\hat{L}=L_{0,0,0}$,  which is realized on the pseudo-euclidean affine 
space.  
The connected component of the group of pseudo-euclidean 
motions is the group topologically generated by the left translations of $\hat{L}$. The 
elements of $\hat{L}$ are the planes 
on which the euclidean metric is induced.        
Two loops  $L_{a,b,c}$ and $L_{a',b',c'}$ in the class ${\cal C}_2$ are 
isomorphic if and 
only if the angles between the tangent $3$-space ${\bf m}_{a,b,c}$, 
respectively   
${\bf m}_{a',b',c'}$ and the tangent $3$-space $\hat{{\bf m}}=
{\bf m}_{0,0,0}$ are the same 
with respect to the Cartan-Killing form. 
This is the case if and only is $a= \pm a'$.
\newline
The non-simple $3$-dimensional differentiable proper Bol loops are either the 
direct products 
of a $1$-dimensional Lie group with a $2$-dimensional Bol loop isotopic to 
the hyperbolic plane loop or the unique Scheerer extension of the Lie group 
$SO_2(\mathbb R)$ by the $2$-dimensional hyperbolic plane loop.   
\end{em}
\newline
\newline
Another impotant class of differentiable loops which generalize the class of  abelian groups is the class of left A-loops. These loops correspond strongly to reductive 
homogeneous spaces. We show that  
any $3$-dimensional differentiable Bol loop which is also a left A-loop is either a Bruck loop or the unique Scheerer extension of the orthogonal group 
$SO_2(\mathbb R)$ by the $2$-dimensional hyperbolic plane loop (which is not a Bruck loop).   
\newline
At the end of this paper we realize within the $3$-dimensional pseudo-euclidean geometry on the 
set of euclidean planes differentiable 
loops which are not Bol loops but the corresponding section of which are given 
in a pure geometric manner.

\section{Bol loops and their Bol algebras}

A set $L$ with a binary operation $(x,y) \mapsto x \cdot y$ is called a loop 
if there exists an element $e \in L$ such that $x=e \cdot x=x \cdot e$ holds 
for all $x \in L$ and the equations $a \cdot y=b$ and $x \cdot a=b$ have 
precisely one solution which we denote by $y=a \backslash b$ and $x=b/a$. 
The left 
translation $\lambda _a: y \mapsto a \cdot y :L \to L$ is a bijection of $L$ 
for any $a \in L$.
Two loops $(L_1, \circ )$ and $(L_2, \ast )$ are called isotopic if there are 
three bijections $\alpha ,\beta ,\gamma : L_1 \to L_2 $ such that 
$\alpha (x) \ast \beta (y)=\gamma (x \circ y)$ holds for any $x,y \in L_1$. 
An isotopism is an equivalence relation. 
If $\alpha =\beta =\gamma $ then the isotopic loops $(L_1, \circ )$ and 
$(L_2, \ast )$ are called isomorphic.
Let $(L_1, \cdot )$ and 
$(L_2, \ast )$ be two loops. The direct product $L=L_1 \times L_2= \{ (a,b)
\ |a \in L_1, b \in L_2 \}$ with the multiplication 
$(a_1,b_1) \circ (a_2,b_2)=(a_1 \cdot a_2, b_1 \ast b_2)$ is again a loop, 
which is called the direct product of $L_1$ and $L_2$, and the loops 
$(L_1, \cdot )$,  $(L_2, \ast )$ are subloops of $(L, \circ )$.   
\newline
 A loop $L$ is called a Bol loop if for any two left translations 
$\lambda _a, \lambda _b$ the product $\lambda _a \lambda _b \lambda _a $ is 
again a left translation of $L$. 
\newline
Let $G$ be the group generated by the left 
translations of $L$ and let $H$ be the stabilizer of $e \in L$ in the group 
$G$. 
The left translations of $L$ form a subset of $G$ acting on the cosets 
$\{x H; x \in G\}$ such that for any given cosets $aH, bH$ there exists 
precisely one left translation $ \lambda _z$ with $ \lambda _z a H=b H$. 
\newline
Conversely let $G$ be a group, H be a subgroup containing no normal 
non-trivial subgroup of $G$ and $\sigma : G/H \to G$ be a section such that 
$\sigma $ satisfies the following conditions:
\newline 
1. The image $\sigma (G/H)$ forms a subset of $G$ with $\sigma (H)=1 \in G$. 
\newline 
2. $\sigma (G/H)$ generates $G$.
\newline
3. $\sigma (G/H)$ acts sharply transitively on the space $G/H$ of the left 
cosets $\{x H, x \in G\}$ (cf. \cite{loops}, p. 18). 
\newline
Then the multiplication on the factor space $G/H$ defined by 
$x H \ast y H=\sigma (x H) y H$ yields a loop $L(\sigma )$. 
This loop is a Bol loop if and only if for all $r,s \in \sigma (G/H)$ the 
element $rsr$ is 
contained in $\sigma (G/H)$.
\newline
If $L_1$ and $L_2$ are Bol loops, then the 
direct product $L_1 \times L_2$ is again a Bol loop.
\begin{Prop}
Let $L$ be a loop
and let $G$ be the group generated by the left translations of $L$, and 
denote by
$H$ the stabilizer of $e \in L$ in $G$. If $G$ and $H$ are direct products 
$G=G_1 \times G_2$ and $H=H_1 \times H_2$ with $H_i \subset G_i$ $(i=1,2)$ 
then the loop $L$ is the direct product of
 two loops $L_1$ and $L_2$,  
and $L_i$ is isomorphic to a loop $L_i^{ \ast }$ having $G_i$ as the group 
generated by its left translations  and $H_i$ as the 
corresponding stabilizer subgroup $(i=1,2)$. 
\newline
In particular there is no $3$-dimensional Bol loop $L$ such that $L$
is the direct product of a $1$-dimensional and a $2$-dimensional Bol loop and 
$G$ is a $5$- or $6$-dimensional Lie group. 
\end{Prop}
\Pro The first assertion is the Proposition 1.18. in \cite{loops}.
The second assertion follows from the fact that a $1$-dimensional Bol loop is a Lie group. Hence   
the dimension of the group $G$ is at most $4$.  
\qed

\bigskip
\noindent
If the elements of $L$ are points of a differentiable manifold and the 
operations $(x,y) \mapsto x \cdot y, \ (x,y) \mapsto x/y, \ (x,y) \mapsto x 
\backslash y
 :L \times L \to L$ are differentiable mappings then $L$ is called a 
differentiable loop. Moreover the manifold $L$ is parallelizable since the set of the left translations is sharply transitive.  
For a differentiable manifold $L$ the group of all 
autohomeomorphisms of $L$ becomes with respect to the compact-open topology  a 
topological group. If $L$ is a Bol loop then the group $G$ topologically 
generated by the left translations within the group of autohomeomorphisms is
 a Lie group (cf. \cite{loops}, p. 33; \cite{quasigroups}, pp. 414-416). 
\newline
Every connected differentiable Bol loop is isomorphic to a Bol loop $L$ 
realized 
on the factor space $G/H$, where $G$ is a connected Lie group, $H$ is a 
connected closed subgroup containing no normal subgroup $\neq \{1\}$ of $G$ 
and the sharply transitive section $\sigma $ of $L$ is a differentiable 
map (cf. \cite{loops}, p. 32). 
\begin{Prop}
If $L$ a connected differentiable $3$-dimensional Bol loop, then 
the pairs $(G,H)$ for the factor space $G/H$ must be different from the following cases:
\newline
a) $G=SL_2(\mathbb C)$ or $PSL_2(\mathbb C)$, $H \in \{U_0,\ U_1 \}$ 
respectively $H \in \{ U_0/ \mathbb Z_2,\ U_1/ \mathbb Z_2 \}$ with 
$U_r=\left \{ \left(\begin{array}{cc}
\displaystyle z & (r-1)w \\ -(r+1)\bar{w} & \bar{z}
\end{array} \right); |z|^2+(r^2-1)|w|^2=1 \right \}$. 
\newline
b) $G$ is locally isomorphic to $SO_3(\mathbb R) \times \mathbb R$ and 
$H$ is a $1$-dimensional subgroup of $G$.  
\newline
c) The group $G=SO_3(\mathbb R) \ltimes \mathbb R^3$ is the connected 
component 
of the euclidean motion group and 
$H$ is a semidirect product of a $2$-dimensional translation group $\mathbb R^2$ by a $1$-dimensional rotation group $SO_2(\mathbb R)$.
\newline
d) $G=SO_3(\mathbb R) \times SO_3(\mathbb R)$, 
and $H$ is any $3$-dimensional subgroup of $G$. 
\newline
e) $G=SL_2(\mathbb C)$ or $PSL_2(\mathbb C)$ and $H=W_r$ respectively 
$W_r \mathbb Z_2/ \mathbb Z_2$, where  $W_r=\left \{ \left(\begin{array}{cc}
\displaystyle \exp((ri-1)x) & 0 \\ z & \exp(-(ri-1)x)
\end{array} \right); x \in \mathbb R, z \in \mathbb C \right \}$ for 
$r \in \mathbb R$. 
\end{Prop}
\Pro In the case b) every $1$-dimensional subgroup $H$ of $G$ containing no 
non-trivial normal subgroup of $G$ has one of the following shapes:
\[ H_1=\{ K \times \{ 0 \} \}, \quad \hbox{or} \quad 
H_2=\{ K \times \varphi (K) \}, \] where $K$ is isomorphic to 
$SO_2(\mathbb R)$ and $\varphi $ is a non-trivial homomorphism. 
\newline
In the case d) every $3$-dimensional subgroup  $H$ of $G$, which does not 
contain any normal subgroup $\neq \{ 1 \}$ of $G$ is conjugate to 
$\{ (a,a) | a \in SO_3(\mathbb R) \}$. 
\newline
In the cases a), c) and for $H_1$ in the case b) the factor space $G/H$ 
is a topological product of spaces having as a factor the $2$-sphere or the projective plane, which are  
 non-parallelizable. In the cases d), e)
and for $H_2$ in the case b) the factor space $G/H$ 
is compact (\cite{loops}, Section 
16). 
\qed

\bigskip
\noindent 
A real vector space $V$ with a trilinear multiplication
$(.,.,.)$ is called a 
Lie triple system $\mathcal{V}$, if the following identities are satisfied:
\newline
\newline
\centerline{ (X,X,Y) =0 } 
\centerline{ (X,Y,Z)+(Y,Z,X)+(Z,X,Y)=0 } 
\centerline{  \big( X,Y,(U,V,W) \big)=\big( (X,Y,U),V,W \big)+\big(U,(X,Y,V),W \big) +\big( U,V,(X,Y,W) \big).}
\newline
\newline
A Bol algebra $A$ is a Lie triple system $(V,(.,.,.))$ with a bilinear 
skew-symmetric operation $[[.,.]]$, $(X,Y) \mapsto [[X,Y]]:V \times V \to V$
such that the following identity is satisfied: 
 \[  [[(X,Y,Z),W]]-[[(X,Y,W),Z]]+(Z,W,[[X,Y]]) \]
 \[ -(X,Y,[[Z,W]])+\big[ \big[ [[X,Y]],[[Z,W]] \big ] \big]=0.  \]       
With any connected differentiable Bol loop $L$ we can associate a Bol algebra 
in the following way: Let $G$ be the Lie group topologically generated by the
 left translations of $L$, and let $(\bf{g, [.,.]}) $ be the Lie algebra of 
$G$. Denote by $\bf{h}$ the Lie algebra of the stabilizer $H$ of $e \in L$ in 
$G$ and by $\bf{m} $ $=T_1 \sigma (G/H)$ the tangent space at $1 \in G$ of 
the image of the section $\sigma :G/H \to G$ corresponding to the Bol loop
 $L$. 
Then $\bf{g}=\bf{m} \oplus \bf{h}$, 
$\big[ [\bf{m},\bf{m}], \bf{m} \big] \subseteq \bf{m}$ and ${\bf m}$ 
generates the Lie algebra 
${\bf g}$. The subspace $\bf{m}$ 
is the tangent Bol algebra of the Bol loop $L$ the operations of which are 
defined by $(X,Y,Z) \mapsto 
\big[ [X,Y], Z \big]$, $(X,Y) \mapsto [X,Y]_{\bf{m}}$, where $X,Y,Z$ 
are elements of the subspace $\bf{m}$ and 
\newline
$Z \mapsto Z_{\bf{m} }: \bf{g} 
\to \bf{m}$ is the projection of $\bf{g}$ onto $\bf{m}$ along the subalgebra 
$\bf{h}$. Miheev and Sabinin define for any Bol algebra $A$ the standard 
enveloping Lie algebra of $A$ (cf. \cite{quasigroups}, pp. 428-429).
If $L$ is a connected differentiable Bol loop having a Bol algebra $A$ as 
its tangent Bol algebra, then 
 the Lie algebra $\bf{g}$ of the Lie group $G$ 
topologically generated by the left translations of $L$ is isomorphic to the
standard enveloping Lie algebra of $A$.  
If $A$ is a $n$-dimensional Bol algebra, then the dimension of its standard 
enveloping Lie algebra  is at most $n+\displaystyle \frac{n(n-1)}{2}$ 
(cf.  \cite{quasigroups}, pp. 428-430). 
\newline
From every triple $(({\bf g}, [.,.]),\ {\bf h},\ {\bf m})$ we can construct  
in a canonical way a triple 
$({(\bf g^*, [.,.]^*), h^*, m})$ (cf. \cite{quasigroups}, pp. 424-425, 
and \cite{loops}, Section 6).  
The Lie algebra $({\bf g^*},[.,.]^*)$ is the Lie algebra of the isometry group of a 
symmetric space and the Lie algebra ${\bf g}$ is an epimorphic image of  
${\bf g^*}$. 
Moreover the subspace ${\bf m} \subset ({\bf g},[.,.])$ determines a Lie 
triple system 
$A:=\big( {\bf m},(.,.,.) \big)$, with $(.,.,.)=\big[ [.,.], .\big]$. If the 
intersection of the centre of ${\bf g}$ with $[{\bf m}, {\bf m}]$ is trivial 
then the Lie algebra ${\bf g}$ is an epimorphic image of the Lie 
algebra of 
the isometry group for the symmetric space belonging to the Lie triple 
system $A$. Hence if ${\bf g}={\bf g}^*$ then ${\bf g}={\bf m} \oplus [ {\bf m}, {\bf m}]$.

\bigskip
\noindent
{\Rem} Let $G$ be a connected Lie group and let ${\bf g}$ be its Lie algebra. For the 
classification of differentiable local Bol loops $L$ 
having $G$ as the group topologically generated by the left translations such that 
${\bf g}= {\bf m} \oplus [{\bf m}, {\bf m}] $ 
for the subspace ${\bf m}$ corresponding to the tangential space of $L$ holds we may proceed in the 
following way: First we determine all symmetric space ${\bf m}$ having $G$ as the isometry group. 
After this we have to find for any symmetric space ${\bf m}$ all subalgebras ${\bf h}$ with ${\bf g}={\bf m} \oplus {\bf h}$ 
such that ${\bf h}$ does not contain any non-trivial ideal of ${\bf g}$. For the classification up to isotopisms it is enough to 
take for ${\bf h}$ a suitable representative of the set $\{ g^{-1} {\bf h} g; g \in G \}$.

\bigskip
\noindent
Since the Bol loops are strongly left 
alternative ( Definition 5.3. in \cite{loops}) every global Bol loop contains an 
exponential image of a 
complement ${\bf m}$ of the Lie algebra ${\bf h}$ of $H$ in the Lie algebra 
${\bf g}$ of $G$, such that ${\bf m}$ generates ${\bf g}$ and satisfies the 
relation $\big [ [{\bf m},{\bf m}], {\bf m} \big ] \subseteq {\bf m}$. 
\newline
\begin{Lemma} Let $L$ be a differentiable loop and denote by ${\bf m}$ the 
tangent space of $L$ at $1 \in G$. Then ${\bf m}$ does not contain any 
 element of $Ad_g {\bf h}$ for some $g \in G$.  
\end{Lemma}
\Pro For $g \in G$ the group $H^g=g^{-1} H g$ is the stabilizer of the
 element $g^{-1}(e) \in L$ in $G$.
\qed

\bigskip
\noindent
The conjugation 
$\sigma (G/H)^g =g ^{-1} \sigma (G/H) g$ of the section $ \sigma $ of $L$
with an element $g \in G$  defines again a sharply transitive 
section $x H \mapsto \sigma (x H)^g$ and the loop $L^*$ belonging to the 
triple $(G,\ H,\ \sigma (G/H)^g )$ 
is isotopic to L (cf. \cite{loops}, Theorem 1.11 (iii)  p. 22).
Moreover this loop $L^*$ is 
isomorphic to $L=(G,\ H,\ \sigma )$ if and only if there is an automorphism 
of $G$ leaving $H$ invariant and mapping the section $\sigma $ corresponding to 
 $L$ onto the section $\sigma (G/H)^g $ corresponding to $L^*$ (cf. 
\cite{loops}, Theorem 1.11 (i)  p. 21).

\bigskip
\noindent
A loop $L$ is called a left A-loop if each $\lambda _{x,y}= \lambda_{xy}^{-1}
\lambda_x \lambda_y:L \to L$ is an automorphism of $L$. The group $G$ 
topologically generated by the left translations of $L$ is a Lie group 
(cf. \cite{loops}, Proposition 5.20, p. 75).
Let ${\bf g}$ be the Lie algebra of $G$ and let ${\bf h}$ be the Lie 
subalgebra of the stabilizer $H$ of $e \in L$ in $G$. Denote by $\sigma :x 
\mapsto \lambda_x :L \to G$ the section corresponding to $L$. Then the 
tangent space ${\bf m}=T_1 \sigma (L)$ of the image of the section 
$\sigma $ at 
$1 \in G$ is complementary to ${\bf h}$ and satisfies the properties  
$ {\bf m} \oplus {\bf h}= {\bf g}$ and $[{\bf h},{\bf m}] 
\subseteq {\bf m}$
(cf. \cite{loops}, Definition 5.18. and Proposition 5.20.  pp. 74-75),  
i.e. ${\bf g}$ is reductive with respect to $({\bf h},{\bf m})$ in the sence 
of Kobayashi and Nomizu (cf. \cite{foundations1} Vol II, p. 190).
\newline
A  differentiable loop $L$ is called a Bruck loop if 
there is an involutory automorphism $\sigma $ of the Lie algebra ${\bf g}$ 
of the connected Lie group 
$G$ generated by the left translations of $L$ such that the tangent space 
$T_e(L)={\bf m}$ is the $-1$-eigenspace and the Lie algebra ${\bf h}$ of the 
stabilizer $H$ of $e \in L$ in $G$ is the $+1$-eigenspace of $\sigma $.

\section{3-dimensional Lie triple systems}

Let $({\bf m}, (.,.,.))$ be a Lie triple system. Denote by $({\bf g}^*, [.,.])$ 
the Lie algebra of the group of displacements of the symmetric space belonging to a 
Lie triple system. We classify the $3$-dimensional non-solvable Lie triple systems 
$({\bf m}, (.,.,.))$ within the group ${\bf g}^*$ as subspaces with $(.,.,.)=[[.,.],.]$. 
\newline
\newline
a) The semisimple Lie triple systems
\newline
If ${\bf m}$ is a $3$-dimensional semisimple Lie triple system then the corresponding Lie 
algebra is 
 semisimple with $4 \le \hbox{dim}\  
{\bf g}^* \le 6$ (\cite{structure}, p. 219, and Theorem 2.7. p. 222). 
Every semisimple Lie algebra ${\bf g}^*$ has dimension $6$. 
To classify the $3$-dimensional semisimple Lie triple systems we have to determine all 
involutory automorphisms of the Lie algebra ${\bf g}^*$, which leave a $3$-dimensional 
subalgebra of ${\bf g}^*$ elementwise fixed. This problem is equivalent to give all involutory 
automorphisms of the Lie group $G^*$, which fix elementwise a $3$-dimensional subgroup of $G^*$.    
\newline
Any $6$-dimensional
 semisimple Lie group is locally 
isomorphic to one of the following Lie groups: 
\newline
1. $PSL_2(\mathbb R ) \times SO_3(\mathbb R ) $
\newline
2. $PSL_2(\mathbb R ) \times PSL_2(\mathbb R )$
\newline
3. $SO_3(\mathbb R ) \times SO_3(\mathbb R )$
\newline
4. $PSL_2(\mathbb C )$.
\newline
{\bf 1. case}: The automorphism group $\Gamma $ of the Lie group $G^*=
PSL_2(\mathbb R ) \times SO_3(\mathbb R ) $ is the direct product of the 
automorphism group of $PSL_2(\mathbb R)$ and the automorphism group of $SO_3(\mathbb R)$. Since 
there is no involutory automorphism of $\Gamma $ leaving a $3$-dimensional subgroup of $G^* $ 
elementwise fixed therefore there is no 
$3$-dimensional Lie triple system corresponding to the Lie group $G^*=
PSL_2(\mathbb R ) \times SO_3(\mathbb R ) $. 
\newline
{\bf 2. case}: The automorphism group $\Gamma $ of the Lie group $G^*=PSL_2(\mathbb R ) \times PSL_2(\mathbb R )$ is the semidirect product 
of the normal automorphism group 
\newline
$\Gamma _1 \times \Gamma _1$, where $\Gamma _1$ is the automorphism group of $PSL_2(\mathbb R)$, 
by the group generated by the 
automorphism $\sigma :G^* \to G^*; \ (u,v) \mapsto (v,u)$. 
The only involutory automorphisms centralizing a $3$-dimensional subgroup of $\Gamma $ are the 
conjugate elements to $\sigma $. Therefore there is 
up to conjugation precisely one $3$-dimensional Lie triple system ${\bf m}$ having the shape 
$\{(X,-X)\  |\  X \in sl_2(\mathbb R)\}$.  
\newline
{\bf 3. case}: In the case $G^*=SO_3(\mathbb R) \times SO_3(\mathbb R)$ there is up to 
conjugation only one Lie triple system corresponding to $G^*$. This Lie triple system has the 
form $\{ (X, -X)\  |\  X \in so_3(\mathbb R) \}$. 
\newline
{\bf 4. case}: {\bf B1} A real basis of the Lie algebra $sl_2(\mathbb C)$: Let us 
consider in ${\bf g}^*$ the following real basis $\{ \cal{H}, \cal{T}, \cal{U}, \hbox{i} 
\cal{H}, \hbox{i} \cal{T}, \hbox{i} \cal{U} \}$, where 
$\cal{ H}=\left( \begin{array}{rr}
1 & 0 \\
0 & -1 \end{array} \right)$, $\cal{ T}=\left( \begin{array}{rr}
0 & 1 \\
1 & 0 \end{array}  \right)$, $\cal{ U}=\left( \begin{array}{rr}
0 & 1 \\
-1 & 0 \end{array}  \right)$. The Lie algebra multiplication is given by
the following rules:
$[\cal{H}, \cal{T}]= \hbox{2} \cal{U},\ [\cal{H}, \cal{U}]=\hbox{2} \cal{T}, 
\ [\cal{U}, \cal{T}]= \hbox{2} \cal{H}$.
\newline
{\bf K1}
The normalized complex Cartan-Killing form $k_{\mathbb C}:sl_2(\mathbb C) 
\times sl_2(\mathbb C) \to \mathbb C$ of $sl_2(\mathbb C)$  
 is the 
bilinear form
defined by $k_{\mathbb C}(X,Y)=\displaystyle \frac {1}{8} 
\hbox{trace} (\hbox{ad}X\ \hbox{ad}Y)$. If $X \in sl_2(\mathbb C)$ 
has the decomposition 
\newline
\centerline{
$X=\lambda _1\  { \cal H}+\lambda_2\ {\cal T}+
\lambda_3\ {\cal U}+ \lambda_4\  \hbox{i} {\cal H}+ \lambda_5\  \hbox{i} 
{\cal T}+\lambda_6\  \hbox{i} {\cal U}$}
 then the complex Cartan-Killing form
$k_{\mathbb C}$ satisfies 
\newline
\centerline{
$k_{\mathbb C}(X)=\lambda _1^2+\lambda_2^2+
\lambda_6^2-
\lambda_3^2- \lambda_4^2- \lambda_5^2+i\ (2\ \lambda _1  \lambda_4+2\ 
\lambda_2 \lambda_5-2\ \lambda_3 \lambda_6)$} 
(cf. \cite{tits1}, Section 6.1, pp. 215-228).
The normalized real Cartan-Killing form $k_{\mathbb R}:sl_2(\mathbb C) \times 
sl_2(\mathbb C) \to \mathbb R$ is the restriction of $k_{\mathbb C}$ to 
$\mathbb R$ such that $k_{\mathbb R}(X)=\lambda _1^2+\lambda_2^2+\lambda_6^2-
\lambda_3^2- \lambda_4^2- \lambda_5^2 $ and the basis $\{ {\cal H}, {\cal T},
{\cal U}, \hbox{i} {\cal H},  \hbox{i} {\cal T},  \hbox{i} {\cal U} \}$ is 
orthonormal with respect to $k_{\mathbb R}$.
\newline
The automorphism group $\Gamma $ of $PSL_2(\mathbb C)$ is the semidirect product of 
the group of inner automorphisms by the group of order 2 generated by the  outer automorphism 
$\tau : z \mapsto \bar{z}$. According to 
(\cite{lowen}, pp. 153-154) there are in $\Gamma $ precisely two conjugacy classes of involutory 
automorphisms centralizing $3$-dimensional subgroups. In the one class there is involutory 
automorphism fixing the group $\left \{ \left ( 
\begin{array}{rr}
a & b \\ 
-\bar{b} & \bar{a}
\end{array} \right); a,b \in \mathbb C, a\bar{a}+b\bar{b}=1 \right \} / 
\mathbb Z_2 \cong SO_3(\mathbb R)$ elementwise. In the other conjugacy class there is an 
involutory automorphism centralizing the subgroup  
$\left \{ \left ( 
\begin{array}{cc}
a & b \\ 
\bar{b} & \bar{a}
\end{array} \right); a,b \in \mathbb C, a\bar{a}-b\bar{b}=1 \right \} / 
\mathbb Z_2 \cong PSL_2(\mathbb R)$. 
\newline 
{\bf 4.1} The subalgebra ${\bf h}_1$ isomorphic to $so_3(\mathbb R)$ is generated by the basis 
elements $\hbox{i} {\cal H}, \hbox{i} {\cal T}, {\cal U}$. The Lie triple system ${\bf m}$ 
with the property $[{\bf m},{\bf m}]= {\bf h}_1$ has as generators  
${\cal H}, {\cal T}, \hbox{i} {\cal U} $. 
\newline
{\bf 4.2} The subalgebra ${\bf h}_2$ isomorphic to $sl_2(\mathbb R)$ has as basis elements  
$\hbox{i} {\cal H}, {\cal T}, \hbox{i} {\cal U}$. The Lie triple system ${\bf m}$ with the 
property $[{\bf m},{\bf m}]= {\bf h}_2$ is generated by the elements 
${\cal H}, \hbox{i} {\cal T}, {\cal U}$. 
\newline
\newline
b) Since there is no non-solvable 
$3$-dimensional Lie triple system with a $2$-dimensional radical 
${\bf r(m)}$ (cf. \cite{structure}, Theorem 2.20. p. 227) we have to consider  
Lie triple systems having a $1$-dimensional radical.  
\newline
Let the subspace $\mathbb R e_1$ be the $1$-dimensional radical ${\bf r(m)}$ of the Lie triple system ${\bf m}$. 
One 
has ${\bf m}={\bf r(m)} \oplus {\bf u_0}$, where ${\bf u_0}$ is a 
$2$-dimensional semisimple subsystem generated by $e_2$ and $e_3$. Since the isometry 
group of the symmetric space ${\cal P}$ corresponding to $u_0$ is $3$-dimensional 
we may assume that $[e_2,e_3]=e_4$ is the generator of the Lie algebra of the  
stabilizer of a point in ${\cal P}$.  
\newline
According to Furness ( \cite{locally}, pp. 44-45) we have $3$ types of  
$2$-dimensional semisimple Lie triple systems. The elements $e_2, e_3$ generate the Lie 
algebra of 
$SO_3(\mathbb R)$ if and only if 
\[ [e_4,e_2]=e_3, \quad [e_4,e_3]=-e_2. \] 
The elements 
 $e_2, e_3$ generate the Lie algebra of 
$PSL_2(\mathbb R)$ if and only if either 
\[ [e_4,e_2]=-e_3, \quad [e_4,e_3]=e_2 \] 
and $e_4$ is elliptic with respect to the Cartan-Killing form, or
\[ [e_4,e_2]=-e_3, \quad [e_4,e_3]=-e_2 \] 
and the element $e_4$ is hyperbolic. 
\newline
Since ${\bf m}$ is a Lie triple system with the properties 
 $\big [ [{\bf m}, e_1], e_1 \big ]=0$ and $\big [ [e_1, {\bf m}], {\bf m} \big ] \subset 
\langle e_1 \rangle $ we 
obtain that
\[1. \quad  [ [ e_1, e_2], e_2]=a e_1,\  [ [ e_1, e_3], e_3]=b e_1,\  
[ [e_1, e_2], e_3]=c e_1 \]
\[  [ [ e_1, e_3], e_2]=d e_1,\   [ [ e_2, e_3], e_1]=(d-c) e_1, \] 
where $a,b,c,d$ are real parameters. 
\newline
{\bf 5. } First let $\mathbb R e_1$ be the centre of ${\bf m}$. Then we have 
$\big [ [e_1, {\bf m}], {\bf m} \big ]=0$ and $a=b=c=d=0$. 
We obtain $3$ types of Lie triple systems with $1$-dimensional 
centre belonging to the different types of the $2$-dimensional semisimple Lie triple systems. Now we characterize the corresponding Lie algebras ${\bf g}^*$. 
\newline
{\bf 5.1} The Lie algebra ${\bf g}^*$ is isomorphic to the Lie algebra 
$so_3(\mathbb R) \oplus \mathbb R$. The multiplication table in ${\bf g}^*$ is defined by the following relations:    
\[ [e_4,e_2]=e_3, \quad [e_4,e_3]=-e_2, \quad [e_4,e_i]=0 \quad (i=2,3,4). \]
\newline
There are precisely two Lie triple systems ${\bf m}_1$, ${\bf m}_2$, such that the corresponding Lie algebra ${\bf g}^*$ is isomorphic to the Lie algebra $sl_2(\mathbb R) \oplus \mathbb R$. 
\newline  
{\bf B2} A real basis of the Lie algebra $sl_2(\mathbb R) \oplus \mathbb R$:   Denote by $e_1$ 
the generator of the Lie algebra $\mathbb R$. Let 
$e_2= ({\cal H}, 0 )$, $e_3= ({\cal T}, 0 )$, $e_4=( {\cal U},0 )$,   
 with ${\cal H}, {\cal T}, {\cal U}$ defined in {\bf B1},  
be a real basis of $sl_2(\mathbb R) \oplus \{ 0 \}$. Then the multiplication in ${\bf g}^*$ 
is given by the following rules: 
\[ [e_2,e_3]=e_4, \ [e_4,e_2]=-e_3, \ [e_4,e_3]=e_2, \ [e_1,e_i]=0 \  (i=2,3,4). \]
{\bf 5.2} The Lie triple system ${\bf m}_1$ is generated by the basis elements 
$e_1, e_2, e_4$. 
\newline
{\bf 5.3} The Lie triple system ${\bf m}_2$ has as generators $e_1, e_2, e_3$. 
\newline
\newline
If the centre of ${\bf m}$ is trivial then not all parameters $a,b,c,d$  
in the system 1) of b) are zero.  
Using the third property in the definition of a 
Lie triple system and an automorphism $\varphi $ of the form: 
$\varphi (e_1)= \alpha \ e_1$, $\varphi (e_2)= \beta \ e_2+ \gamma \ e_3$,  
$\varphi (e_3)=\delta \ e_2+\varepsilon e_3$, where 
$ \beta \varepsilon - \gamma \delta \neq 0$, we have the following four cases: 
\newline
{\bf 6.1} The Lie triple system ${\bf m}_1$ is given as follows:  
\[ [e_4, e_2]=e_3,\ [ e_4, e_3]=-e_2,\ [[e_1, e_2], e_2]=-e_1,\ 
[[e_1, e_3], e_3]=-e_1. \]  
The corresponding Lie algebra ${\bf g}^*$ is isomorphic to the Lie algebra $so_3(\mathbb R) \ltimes \mathbb R^3$. Since $
\mathbb R e_2,\ \mathbb R e_3,\ \mathbb R e_4$ correspond to 
$1$-dimensional rotations and $[[e_1,e_3],e_2]=0,\ [[e_1,e_2],e_3]=0,\ 
[e_4,e_1]=0$ the vectors $[e_1,e_3],\ [e_1,e_2],\ e_1$ are the axes of the 
rotation groups corresponding to $e_2$, $e_3$ respectively $e_4$. Hence 
$e_1$, $e_5:=[e_1,e_3]$, $e_6:=[e_1,e_2]$ form a basis of the radical 
of ${\bf g}^* \cong \mathbb R^3$. We have  
 $[e_4,e_6]=e_5$, $[e_4,e_5]=-e_6$ because the rotation group belonging to 
$ e_4$ with the axis $ e_1 $ leaves the 
$2$-dimensional subspace $\langle e_5, e_6 \rangle$ orthogonal to 
$\langle e_1 \rangle$ invariant. Similarly the rotation groups corresponding 
to $e_2$  and to $e_3$ with axes $e_5$ respectively $e_6$ leaves the 
subspaces $\langle e_1, e_6 \rangle$, respectively 
$\langle e_1, e_5 \rangle$ invariant. Hence we have $[e_2, e_6]=e_1$, 
$[e_3,e_5]=e_1$. Therefore the multiplication in the Lie algebra ${\bf g}^*$ is 
determined. 
\newline
{\bf B3} Representing the Lie algebra ${\bf g}^*$ in the Lie algebra of real 
$(4 \times 4)$-matrices we may choose as basis the following matrices: 
\[ e_1=\left ( \begin{array}{rrrr}
0 & 0 & -1 & 0 \\
0 & 0 & 0 & 0 \\
0 & 0 & 0 & 0 \\
0 & 0 & 0 & 0 \end{array} \right ),\ e_2= \left ( \begin{array}{rrrr}
0 & 0 & 0 & 0 \\
0 & 0 & 0 & 0 \\
0 & 0 & 0 & -1 \\
0 & 0 & 1 & 0 \end{array} \right ),\  e_3=\left ( \begin{array}{rrrr}
0 & 0 & 0 & 0 \\
0 & 0 & 1 & 0 \\
0 & -1 & 0 & 0 \\
0 & 0 & 0 & 0 \end{array} \right ), \]
\[  \  e_4=\left ( \begin{array}{rrrr}
0 & 0 & 0 & 0 \\
0 & 0 & 0 & 1 \\
0 & 0 & 0 & 0 \\
0 & -1 & 0 & 0 \end{array} \right ),
\  e_5=\left ( \begin{array}{cccc}
0 & 1 & 0 & 0 \\
0 & 0 & 0 & 0 \\
0 & 0 & 0 & 0 \\
0 & 0 & 0 & 0 \end{array} \right ), 
e_6=\left ( \begin{array}{cccc}
0 & 0 & 0 & 1 \\
0 & 0 & 0 & 0 \\
0 & 0 & 0 & 0 \\
0 & 0 & 0 & 0 \end{array} \right ). \]
{\bf 6.2} 
\[  [ e_4, e_2]=-e_3,\ [ e_4, e_3]=e_2,\ [[e_1, e_2], e_2]=e_1,\ 
[[e_1, e_3], e_3]=e_1. \]
{\bf 6.3} The Lie triple ${\bf m}_3$ is given by the following relations: 
\[  [e_4, e_2]=-e_3,\ [e_4, e_3]=-e_2,
\ [[e_1, e_2], e_2]=e_1,\ [[e_1, e_3], e_3]=-e_1. \]
In both cases the Lie algebra ${\bf g}^*$ is isomorphic to $sl_2(\mathbb R) \ltimes \mathbb R^3$. From the relation in {\bf 6.3} we can deduce the multiplication table of the Lie algebra ${\bf g}^*$:  
\[ [e_1,e_2]=:e_6,\ [e_1,e_3]=:e_5,\ [e_2,e_3]=:e_4,\ [e_5,e_4]=-e_6, \]
\[ [e_1,e_4]= [e_1,e_5]=[e_1,e_6]= [e_2,e_5]= [e_3, e_6]=[e_6,e_5]=0,\] 
\[ [e_2,e_6]= [e_3,e_5]= -e_1, \ [e_2,e_4]=e_3,\ [e_3,e_4]=-e_2,\  
[e_6,e_4]=e_5. \] 
{\bf B4} Representing ${\bf g}^*$ as pairs of real $(2 \times 2)$-matrices we may choose as basis 
$e_1=\left( 0, -{\cal U} \right)$, $e_2=\left( {\cal H},0 \right)$, 
$e_3=\left( {\cal T}, 0 \right)$, $e_4=\left( {\cal U}, 0 \right)$,  
$e_5=\left( 0, -{\cal H} \right)$, $e_6=\left( 0, {\cal T} \right)$ with 
${\cal H}, {\cal T}, {\cal U}$ defined in {\bf B1}.  
The Lie triple system ${\bf m}_2$ defined by the relations {\bf 6.2} in the basis {\bf B4} is 
generated by the basis elements $e_2,e_4,e_6$. 
The Lie triple system ${\bf m}_3$ defined by the relations {\bf 6.3} belonging to the basis 
{\bf B4} has as generators $e_1,e_2,e_3$.
\newline
If an element $X$ of ${\bf g}^*=sl_2(\mathbb R) \ltimes \mathbb R^3$ has the 
decomposition \[ X=\lambda _1\ e_1+ \lambda _2\ e_2+ \lambda_3\ e_3+ 
\lambda_4\ e_4+ \lambda_5\ e_5+ \lambda_6\ e_6 \] then the Cartan-Killing 
form $k$ (cf.  \cite{tits1}, p.193)  on ${\bf g}^*$ satisfies the relation
\[ \ {\bf K2}: \quad \quad k(X,X)=\lambda_2^2+ \lambda_3^2 -\lambda_4^2. \] 
{\bf 7} \quad The last Lie triple system and the corresponding Lie algebra ${\bf g}^*$ are defined 
by the following rules:  
\[ [e_4, e_2]=-e_3,\ [e_4, e_3]=-e_2,\ [[e_1, e_2], e_2]=
\displaystyle \frac{1}{4}\ e_1, 
 [[e_1, e_3], e_3]=- \displaystyle \frac{1}{4}\ e_1,\] \[ [[e_1, e_2], e_3]=
\displaystyle \frac{1}{4}\ e_1,\ 
[[e_1, e_3], e_2]=-\displaystyle \frac{1}{4}\ e_1,\ [[e_2, e_3], e_1]=-
\displaystyle \frac{1}{2}\ e_1. \]
In this case the Lie algebra ${\bf g}^*$ is isomorphic to 
$sl_2(\mathbb R) \ltimes \mathbb R^2$. 
Moreover the elements $e_2, e_3$ and $e_4=[e_2,e_3]$ form a basis of $sl_2(\mathbb R)$, such 
that $e_2$ is a hyperbolic, $e_3$ an elliptic and $e_4$ again a hyperbolic element with repect 
to the Cartan-Killing form.

\section{Bol loops having  semisimple Lie groups as the groups generated by 
their left translations}

In this section we classify all 3-dimensional connected  differentiable Bol 
loops having semisimple Lie groups as the groups topologically generated by their 
left translations and  
describe the symmetric spaces as well as the natural geometries associated with these loops.  
\newline
Let $L$ be a 3-dimensional proper Bol loop, let $G$ be the group generated by 
the left 
translations of $L$ and $H$ be the stabilizer of $e \in L$ in $G$. 
Since $G$ is semisimple and $4 \le  \hbox{dim}\  G \le 6$ we have 
$ \hbox{dim}\  G = 6$.

\begin{Prop} Let $G=G_1 \times G_2$ be 
topologically generated by the left translations of a $3$-dimensional 
connected differentiable Bol loop $L$ such that  $G_i $ 
\newline
$(i=1,2)$
are $3$-dimensional 
quasi-simple 
Lie groups. Then $G_i$ $(i=1,2)$ is 
isomorphic  to $PSL_2(\mathbb R)$ and  the 
stabilizer $H$ of $e \in L$ in $G$ may be chosen either 
as 
\newline
(i) \centerline{$H_1=\{ ( x,x )\ |\  x \in PSL_2(\mathbb R) \}$} 
or 
\newline
(ii) \centerline{$H_2=\left\{ \left( \left( \begin{array}{ll}
a & b_1 \\
0 & a^{-1} \end{array} \right), \left( \begin{array}{ll} 
a & b_2 \\
0 & a^{-1} \end{array} \right) \right); a > 0, b_1,b_2 \in \mathbb R \right 
\}$. }  
\end{Prop}
\Pro  Denote by $\pi_i:G \to G_i$ the natural projection of $G$ to $G_i$ for 
$i=1,2$. Let 
$\hbox{dim}\  \pi_1(H) \le 1$. Since $H \le \pi_1(H) \times \pi_2(H)$ one has 
$\hbox{dim}\  \pi_2(H) \ge 2$. If $\hbox{dim}\  \pi_2(H) = 2$ then $H$ is the 
direct product of $\pi_1(H)$ with $\pi_2(H)$ and we obtain a contradiction to 
Proposition 1. Let $\pi_2(H)=G_2$.   
If $\hbox{dim}\  \pi_1(H)=0$ then 
$H=\{1 \} \times G_2$ which is impossible. If $\pi_1(H)$ is a $1$-dimensional 
subgroup then $H$ would be $4$-dimensional since there is no non-trivial 
homomorphism from $G_2$ into $\pi_1(H)$. 
\newline
Let now $\hbox{dim}\  \pi_1(H) = 2$. 
We may assume that $\hbox{dim}\  \pi_2(H) \ge 2$ since interchanging the indices 
we would obtain the previous case. Therefore each of the factors of $G$ is 
locally isomorphic to the group $PSL_2(\mathbb R)$. Since the Lie algebra of 
$G_2$ is simple  there is no 
non-trivial homomorphism from $\pi_2(H)=G_2$ into $\pi_1(H)$, so $\pi_2(H)$ cannot be 
$G_2$. If 
$\hbox{dim}\  \pi_2(H) = \hbox{dim}\  \pi_1(H)=2$ then 
$\pi_i(H) \cong {\cal L}_2=\{ a x +b\  |\  a>0, b \in \mathbb R \}$ 
$(i=1,2)$ and 
there 
exist homomorphisms $\varphi _1, \varphi _2$ with $1$-dimensional kernels, 
such that
$\varphi _1 :  \pi_1(H) \to \pi_2(H)$ and $\varphi _2 :  \pi_2(H) \to 
\pi_1(H)$. 
Since $H \cap G_i= ker\  \varphi _i=\left \{ \left ( \begin{array}{cc}
1 & b \\
0 & 1 \end{array} \right ), b \in \mathbb R \right \}$ and (up to conjugation) 
$im\ \varphi _i=
\left \{ \left ( \begin{array}{ll}
a & 0 \\
0 & a^{-1} \end{array} \right ), a>0 \right \}$ $(i=1,2)$ we obtain that 
\newline
\centerline{
$\pi_1(H)=ker\ \varphi_1 \ im \ \varphi _2$, 
$\pi_2(H)=ker\ \varphi_2 \ im \ \varphi _1$} and $H$ has the shape (ii). 
\newline
Finally let $\hbox{dim} \pi_1(H)=3$. From the previous  
arguments it follows that 
\newline
$\hbox{dim}\  \pi_2(H)$ 
$=3$. If there is no 
homomorphism $\varphi:G_1 \to G_2$ then there is no 
$3$-dimensional subgroup $H$ of $G=G_1 \times G_2$. If there exists a 
homomorphism $\varphi : G_1 \to G_2$ then the stabilizer $H$ has the shape 
$\{ (x, \varphi (x))\  |\  x \in G_1 \}$. Moreover $\varphi $ must be an isomorphism since 
otherwise $G$ would contain a discrete central 
subgroup of $G$. 
Hence  $G$ can be 
identified with 
$G_1 \times G_1$, where $G_1$ is isomorphic either to $PSL_2(\mathbb R)$ or 
$SO_3(\mathbb R)$, and $H$ may be chosen as the diagonal subgroup 
$\{ (x,x)\  |\  x \in G_1 \}$. According to Proposition 2 d) we have  
$G=PSL_2(\mathbb R) \times PSL_2(\mathbb R)$. 
\qed

\bigskip
\noindent 
Now let $G$ be locally isomorphic to the group $PSL_2 (\mathbb C)$. 
According to ( \cite{Asoh},  pp. 273-278) there are $4$ conjugacy classes of 
the $3$-dimensional subgroups of $G=SL_2(\mathbb C)$, which we denote by 
$W_r$, $U_0$, $U_1$ and $SU_2(\mathbb C)$. 
Since $SU_2(\mathbb C)$ contains central elements $\neq 1$ of $SL_2(\mathbb C)$ it follows 
from Proposition 2 that no of these groups can be the stabilizer of $e \in L$ in $G$. 
Hence it is sufficient to consider the 
following  cases: 
\newline
1. $G \cong PSL_2(\mathbb R ) \times PSL_2(\mathbb R )$, 
and $H_i$ $(i=1,2)$ is one of the subgroups given in Proposition 4. 
\newline
2. $G$ is isomorphic to $PSL_2(\mathbb C )$ and $H$ is  
isomorphic to $SO_3(\mathbb R)$.

\bigskip
\noindent
Now we deal with the case 1). There is up to isomorphism precisely one Lie triple system 
${\bf m}=\{ (X,-X)\  |\  X \in sl_2(\mathbb R) \}$ corresponding to the Lie group 
$PSL_2(\mathbb R) \times PSL_2(\mathbb R)$ ({\bf 2. case}, Section 3). 
\newline
The Lie algebra of $H_1$ is ${\bf h}_1=\{ (X,X)\  |\  X \in sl_2(\mathbb R) \}$. 
The Lie algebra ${\bf h}_2$ of $H_2$ is generated by the basis elements: 
$({\cal H}, {\cal H})$, $({\cal U} + {\cal T},0)$, $(0,{\cal U}+ {\cal T})$, 
where     
${\cal H}, \ {\cal U}, \ {\cal T}$ are $(2 \times 2)$-matrices defined in B1. 
We see that the intersection of ${\bf m}$ and 
${\bf h}_1$ is trivial. But the subalgebra ${\bf h}_1$ 
contains the element $({\cal H}, {\cal H})$, which is conjugate to the element 
$({\cal H},- {\cal H}) \in {\bf m}$ under the element  
corresponding  
 to the pair of matrices $\left( \pm 1, \pm \left( 
\begin{array}{rr}
0 & 1 \\
-1 & 0 \end{array} \right) \right)$. 
This contradicts Lemma 3. 
\newline
Since the intersection of ${\bf m}$ and ${\bf h}_2$ is 
$\langle ({\cal U} + {\cal T},-({\cal U} + {\cal T})) \rangle $ the subspace ${\bf m}$ cannot be a complement 
to ${\bf h}_2$ in ${\bf g}$. 
\newline
From this discussion it follows that there is no $3$-dimensional 
connected 
differentiable Bol loop such that the group $G$ generated by its left 
translations is 
isomorphic to the group $PSL_2 (\mathbb R) \times PSL_2( \mathbb R)$.

\bigskip
\noindent
Now we consider the case 2. 
The Lie algebra ${\bf h}$ of $H$ with respect to the basis {\bf B1}  ({\bf 4.case}, Section 3) 
 has as generators $\cal{U}, \hbox{i} \cal{T}, \hbox{i} \cal{H} $. The Lie triple system 
${\bf m}=\langle \hbox{i} \cal{T}, \cal{U}, \cal{H} \rangle$ given in {\bf 4.2} cannot be a 
complement 
to ${\bf h}$ in ${\bf g}$ since ${\bf h} \cap {\bf m}=\langle  \hbox{i} \cal{T}, \cal{U} \rangle$. 
\newline
The Lie triple system ${\bf m}$ of {\bf 4.1} satisfies $[{\bf m}, {\bf m}]$
$={\bf h}$, 
 and 
$\bf{g}= {\bf m} \oplus [{\bf m},{\bf m}]$. 
Hence it determines a $3$-dimensional Riemann symmetric space (cf. \cite{foundations1}, Chapter 
VI, Theorem 2.2 (iii)). According to Theorem 11.8 in (\cite{loops}) there exists a global 
differentiable Bruck loop $L_0$ homeomorphic to $\mathbb R^3$ and having 
$G \cong PSL_2(\mathbb C)$ 
as the group topologically generated by its left translations such that $H \cong SO_3(\mathbb R)$
is the stabilizer of $e \in L$ in $G$. 
Since the symmetric spaces with $PSL_2(\mathbb C)$ as the group of displacements which correspond 
to global differentiable Bol loops are isomorphic, there is precisely one isotopism class 
${\cal C}$ of differentiable Bol loops having $PSL_2(\mathbb R)$ as the group topologically 
generated by their left translations (cf. Remark in section 3) and $L_0$ is a representative of 
${\cal C}$. Since the subspace ${\bf m}$ corresponding to $L_0$ is the hyperbolic space we call 
$L_0$ the hyperbolic space loop.
\newline
In order to determine within the isotopism class ${\cal C}$ the isomorphism classes we look for a suitable parametrization 
of the local Bol loops isotopic to the hyperbolic space loop. For this reason we consider all $3$-dimensional complements to 
${\bf h}=so_3(\mathbb R)$ in ${\bf g}=sl_2(\mathbb C)$ with the properties 
$\bf{g} = \bf{m} \oplus \bf{h}$, $\big[ [ \bf{m}, \bf{m}], \bf{m} \big] 
\subseteq \bf{m}$ and $\bf{m}$ generates $\bf{g}$.
We can write 
$\bf{m}$ in the general form:
\[ \bf{m}=\langle \cal{T}+\hbox{a} \cal{U}+\hbox{b i} \cal{T}+\hbox{c i} 
\cal{H}, 
\hbox{i} \cal{U}+ \hbox{d} \cal{U}+\hbox{e i} \cal{T}+\hbox{f i} \cal{H},
\cal{H}+ \hbox{g} \cal{U}+\hbox{h i} \cal{T}+ \hbox{k i} \cal{H} \rangle ,  \] 
$a,b,c,d,e,f,g,h,k \in \mathbb R$. According to 
 \cite{sabinin} (pp. 217-219)  there 
are precisely two classes of the Bol complements ${\bf m}$ to $\bf{h}$ in 
$\bf{g}$ having the 
following shapes: 
\newline
${\bf m}_{a}=\langle \cal{T}+ \hbox{a} \cal{U}, \hbox{i} \cal{U}+ 
\hbox{a i} 
\cal{T}, \cal{H} \rangle $ for $ a \in \mathbb R \backslash \{1,-1 \}$, 
\newline
${\bf m}_{d}=\langle \cal{U+T}, \hbox{d i} \cal{H} +\hbox{i} 
(\cal{U+T}), \cal{H}+\hbox{d}\  \cal{U} \rangle $ for $ d \in \mathbb R 
\backslash \{0 \}$.   
\newline
A Bol complement belongs to a local Bol loop in the class ${\cal C}$ 
if and only if $[{\bf m}, {\bf m}]$ is a compact subalgebra of ${\bf g}$. Since the value of the 
complex Cartan-Killing form {\bf K1} on the element 
$[\cal{U+T}, \hbox{d i} \cal{H}+\hbox{i} (\cal{U+T})]$$\in [{\bf m}_{d}, {\bf m}_{d}]$
is zero $[{\bf m}_{d}, {\bf m}_{d}]$ is not compact. The subalgebra $[{\bf m}_{a}, {\bf m}_{a}]$
 is compact if and only if $a \in (-1, 1)$ since 
$[{\bf m}_{a}, {\bf m}_{a}]=\langle \hbox{ i} \cal{H}, \cal{U}+ \hbox{a} \cal{T}, \hbox{i} \cal{T}+ 
\hbox{a i} \cal{U} \rangle$. 
\newline
In the isotopism class ${\cal C}=\{ L_a, a \in (-1, 1) \}$ only the Bol complement ${\bf m}_0$ satisfies the relation $[{\bf h}, {\bf m}_0 ] \subseteq 
{\bf m}_0$. Therefore only the hyperbolic space loop $L_0$ corresponding to the triple 
$(G, H, \exp {\bf m}_0)$ is a left A-loop. It is the unique Bruck loop in the class ${\cal C}$.
\newline
The complement ${\bf m}_0=\langle \cal{T}, \hbox{i} \cal{U}, \cal{H} \rangle $
is orthogonal to $\bf{h}$ with respect to the Cartan-Killing form 
$k_{\mathbb R}$ on
$\bf{g}$ (see {\bf K1} in Section 3). 
Let ${\bf m}_a \subset {\bf g}=sl_2(\mathbb C)$ be the tangent space $T_1 
\sigma_a(G/H)$ of the image $\sigma_a(G/H)$ of the differentiable section 
$\sigma_a:G/H \to G$ belonging to the loop $L_a$ for $a \in (-1,1)$.
Two loops corresponding to $(G,H, \exp {\bf m}_a)$ and $(G,H, \exp {\bf m}_b)$ 
are isomorphic if and only if there exists an automorphism $\alpha $ of 
${\bf g}$ such that $\alpha ({\bf m}_a)={\bf m}_b$ and 
$\alpha ({\bf h})={\bf h}$. The automorphism group of ${\bf g}$ leaving 
${\bf m}_0$ and ${\bf h}$ invariant is the semidirect product $\Theta $ of 
$\hbox{Ad}_H$ by the group generated by the involutory map $\varphi :z \mapsto \bar{z}$. The 
Cartan-Killing 
form $k_{\mathbb R}$ is invariant under $\Theta $. Hence the $3$ angles 
between ${\bf m_a}$ and ${\bf m_0}$ respectively between ${\bf m_b}$ and 
${\bf m_0}$ are equal with respect to $k_{\mathbb R}$. The loops $L_a$ and $L_b$ are isomorphic if there exists 
a matrix $B=\left( \begin{array}{rr}
c_1+c_2\ i & d_1+d_2\ i \\
-d_1+d_2\ i & c_1-c_2\ i \end{array} \right) \in H$ such that 
$\varphi ( B^{-1}\ {\bf m_a}\ B)={\bf m_b}$. This condition is equivalent 
with the following equations:
\[ c_1^2+c_2^2+d_1^2+d_2^2=1,\ \  d_1\ c_2-c_1\ d_2=0, \ \ d_1\ c_1+c_2\ d_2
=0, \] 
\[ c_1\ d_2+c_2\ d_1=0, \ \ c_2\ d_2-c_1\ d_1=0, \]
\[ (d_2^2-d_1^2 -a\ (c_2^2-c_1^2))\ b=c_1^2-c_2^2-a\ (d_2^2-d_1^2), \] 
\[ (d_1\ d_2+c_1\ c_2\ a)\ b=c_1\ c_2-a\ d_1\ d_2, \] 
\[ b\ (-a\ (c_1^2-c_2^2)+(d_2^2-d_1^2))=-a\ (d_2^2-d_1^2)-(c_1^2-c_2^2), \]
\[ b\ (c_1\ c_2\ a-d_1\ d_2)=d_1\ d_2\ a+c_1\ c_2. \]
The solution of this  system of equations is $b=\displaystyle \frac{1}{a}$ and  
$b=-a$. Since $b,a \in (-1,1)$ it follows that $b=-a$.
Therefore the loops $L_a,\ a \in [0,1)$ can be chosen as the representatives 
of the isomorphism classes of the Bol loops $L_a,\ a \in (-1,1)$.

\bigskip
\noindent     
The multiplication of the hyperbolic space loop $L_0$ may be given in the 
$3$-dimensional hyperbolic geometry $\mathbb H _3$ explicitly. We choose a point
$e$ in $\mathbb H _3$ and define 
$x \circ y= \tau _{e,x} (y)$ 
for all points $x,y \in \mathbb H _3$, where $\tau_{e,x}$ is the hyperbolic translation 
along the line joining $e$ with $x$ and mapping $e$ onto $x$.    
\newline
The hyperbolic space geometry has an elementary model in the upper half 
space 
$\mathbb R^{3+}=\{ (x,y,z) \in \mathbb R^3; z>0 \}$. 
We can identificate the elements of $\mathbb R^3$ with the elements of the 
$\bf{J}$-quaternion space. The $\bf{J}$-quaternion space is the 
$3$-dimensional subspace of the quaternion space which is orthogonal to the 
canonical basis quaternion $k$ (cf. \cite{elementary}, p. 3). We can give the
 action of the group $SL_2( \mathbb C)$ on $\mathbb R^3$ by the linear 
rational functions: $\gamma (w)=(a w+b)(c w+d)^{-1}$, where 
$\gamma =\left( \begin{array}{cc}
a & b \\
c & d \end{array} \right)$, $a,b,c,d \in \mathbb C, ad-bc=1, w=x+j y \in 
{\bf J}, x \in \mathbb C$, $y \in \mathbb R$ 
(cf. \cite{elementary}, p. 26).    
The restriction of this action onto the subspace $\mathbb R^{3+}$ defines the 
action of $SL_2(\mathbb C)$ on the upper half space, and $(w,\pm \gamma) 
\mapsto \gamma (w)$ is the transitive action of $PSL_2( \mathbb C)$ on the 
upper half space.
Since the stabilizer subgroup $H=SO_3( \mathbb R)$ leaves 
the point $j$ fixed, we may choose this point as the identity element of the 
hyperbolic space loop. 
Summarizing our discussion we obtain
\begin{Theo} There is only one isotopism class $\cal C$ of the $3$-dimensional
connected differentiable Bol loops $L$ such that the group $G$ generated by 
the left translations $\{ \lambda_x; x \in L \}$ is a semisimple Lie group.
The group $G$ is isomorphic to $PSL_2(\mathbb C)$ and the stabilizer $H$ of 
$e \in L$ in $G$ is isomorphic to $SO_3(\mathbb R)$. 
\newline
Any loop $L_a \quad 
(-1 < a < 1)$ in this class $\cal C$ is isotopic to the hyperbolic space loop
$L_0$. The tangent space $T_1 \Lambda $ of the set $\Lambda $ of the left 
translations of the hyperbolic space loop $L_0$ at the identity $1 \in G$
is the plane ${\bf m_0}$ through $0$ in the Lie algebra ${\bf g}$ of $G$ such 
that ${\bf m}_0$ is orthogonal to the $3$-dimensional Lie algebra ${\bf h}$ 
of 
$H$ with respect to the Cartan-Killing form $k_{\mathbb R}$ of ${\bf g}$. 
\newline
In the class $\cal C$ only the hyperbolic space loop is a Bruck loop.
\newline
The tangential spaces $T_1 \sigma _a (G/H)$ at $1 \in G$ for all 
$a \in (-1, 1)$ are the Bol complements ${\bf m}_a$.
The loops $L_a$ and $L_b$ in  ${\cal C}$ are isomorphic if and 
only if the angles between the $3$-spaces ${\bf m}_a$ and ${\bf m}_0$, 
respectively ${\bf m}_b$ and ${\bf m}_0$ are the same 
with respect to $k_{\mathbb R}$.  
As representatives of the isomorphism classes of the loops belonging to 
 ${\cal C}$ we may choose the loops $L_a$, where 
$0 \le a <1$. 
\end{Theo}

\section{3-dimensional Bol loops corresponding to 
\newline 4-dimensional non-solvable Lie groups }

In this section we determine all 
 $3$-dimensional connected differentiable global Bol loops 
$L$ having a $4$-dimensional non-solvable Lie group $G$ as the group 
topologically generated by their left translations. 
\newline
Since according to Proposition 2 b) a group locally isomorphic 
to $SO_3(\mathbb R) \ltimes \mathbb R$ cannot be  
the group topologically generated by the left translations of a $3$-dimensional Bol loop
we have to consider only the case that $G$ is locally isomorphic to 
$PSL_2(\mathbb R) \times \mathbb R$. 
Moreover, since the group $G$ is isomorphic to the Lie group $G^*$ described in  {\bf 5.2} 
and {\bf 5.3} of section 3 we have to discuss the following situation. 
\newline
Using the $\mathbb R$-basis {\bf B2} ({\bf 5.2}, Section 3) the $1$-dimensional subalgebras ${\bf h}$ of ${\bf g}=sl_2(\mathbb R) \times \mathbb R$ 
containing no non-zero ideal of ${\bf g}$ are the following: 
\newline
${\bf h}_1=\langle e_2+k e_1 \rangle$ 
\newline
${\bf h}_2=\langle (e_3+e_4)+k e_1 \rangle$
\newline
${\bf h}_3=\langle e_4+k e_1 \rangle$, where $k \in \mathbb R$. 
\newline
Using the automorphism $A: {\bf g} \to {\bf g}$ given by $A(e_1)= k e_1$, $A(e_i)=e_i$ for 
$2 \le i \le 4$   we may assume that 
$k=0,1$. 
\newline
Since the Lie triple system ${\bf m}_1$ in the case {\bf 5.2} is generated by the basis 
elements $e_1, e_2, e_4$, the Lie algebras ${\bf h}_1$ and ${\bf h}_3$ cannot be complements 
to ${\bf m}_1$ in ${\bf g}$. Since the elements $e_2+e_4 \in {\bf m}$ and $e_3+e_4 \in {\bf h}$ 
are both parabolic in $sl_2(\mathbb R)$, we have a contradiction to Lemma 3. 
\newline
Now we consider the Lie triple system ${\bf m}_2=\langle e_1,e_2,e_3 \rangle$ in the case 
{\bf 5.3}. The intersection of ${\bf h}_1$ and ${\bf m}_2$ is not trivial, therefore 
${\bf m}_2$ is not a complement to ${\bf h}_1$ in ${\bf g}$.
Since ${\bf h}_2 \cap {\bf m}_2={\bf h}_3 \cap {\bf m}_2=\{ 0 \}$ and there is no element of 
${\bf m}_2$ which is conjugate to an element of ${\bf h}_i$ for $i=2,3$ we have to deal with 
both cases. 
\newline
For $k=0$ we have $H_2=\exp {\bf h}_2 = 
\left \{ \left ( \left ( \begin{array}{cc} 
1 & b \\
0 & 1 \end{array} \right ), 1 \right), b \in \mathbb R \right \}$ and 
\newlineû
$H_3=\exp {\bf h}_3 = \{ (x,1) , x \in SO_2(\mathbb R) \}$. Now  
$G$ and $H$ are direct products $G=G_1 \times G_2$, $H=H_1 \times H_2$ with $H_i \subset G_i$ 
$(i=1,2)$. The Bol loop $L$ corresponding to a such pair is the direct product of a 
$2$-dimensional Bol loop $L_1$ and a $1$-dimensional Lie group $L_2$ such that the group 
topologically generated by the left translations of the $2$-dimensional Bol loop $L_1$ is 
 isomorphic to 
$PSL_2(\mathbb R)$. Then $G=PSL_2(\mathbb R) \times G_2$, where $G_2$ is either  
$\mathbb R$ or $SO_2(\mathbb R)$, and $H=H_3$. Moreover $L_1$ is isotopic to the hyperbolic 
plane loop.  
\newline
For $k=1$ one has $G=PSL_2(\mathbb R) \times \mathbb R$, $H_2=\exp {\bf h}_2 = 
\left \{ \left ( \left ( \begin{array}{cc} 
1 & b \\
0 & 1 \end{array} \right ), b \right) \right \}$, where $b \in \mathbb R $ and 
$\exp {\bf m}= \left \{ \left ( \left ( \begin{array}{cc}
x+y  & z \\
z & x-y \end{array} \right ), k \right ) \right \}$ with 
$x,y,z,k \in \mathbb R$, and $x^2-y^2-z^2=1$   
(\cite{old}, p. 24). If for a section $\sigma: G/H \to G$ the 
image $\sigma(G/H)$ contains $\exp{\bf m}$ then any coset $\left( \left( 
\begin{array}{cc}
1+c & 1 \\
c   & 1 \end{array} \right),0 \right) H,\  c \in \mathbb R$ should contain 
precisely one element $s$ of $\sigma(G/H)$. If $c \neq -1$ then 
$s=\left( \left( \begin{array}{cc}
1+c & c \\
c   & \displaystyle \frac{c^2+1}{1+c} \end{array} \right),0 \right)$. Since $\sigma: G/H 
\to G$ is continuous we have 
$ \sigma \left( \left( \left( 
\begin{array}{cc}
0 & 1 \\
-1   & 1 \end{array} \right),0 \right) H \right)=
\lim_{c \to -1} \sigma \left( \left( \left( 
\begin{array}{cc}
1+c & 1 \\
c   & 1 \end{array} \right),0 \right) H \right)= $
\newline
$ \lim_{c \to -1} \left( \left( \begin{array}{cc}
1+c & c \\
c   & \displaystyle \frac{c^2+1}{1+c} \end{array} \right),0 \right)$ ,
which is a contradiction. This means that there is no global Bol loop 
corresponding to the triple $(G,H_2, \exp {\bf m})$.
\newline
The other case is $G=PSL_2(\mathbb R) \times SO_2(\mathbb R)$, 
$H_3=\{ (x,x)\ |\ x \in SO_2(\mathbb R) \}$ and 
\newline
$\exp {\bf m} = \exp \{ ( \lambda_2 e_2+ \lambda_3 e_3), \lambda_2 ,\lambda_3 \in \mathbb R \} \times \exp \{ 
(\lambda_1 e_1), \lambda_1 \in \mathbb R \}=M \times SO_2(\mathbb R)$,
such that $M$ is the image of the section  
$\sigma _1$ given by 
\[ \left( \begin{array}{cc}
a & 0 \\
b & a^{-1} \end{array} \right) \mapsto \left ( \begin{array}{cc}
a & 0 \\
b & a^{-1} \end{array} \right) \left ( \begin{array}{cc}
\frac{a^{-1}+a}{\pm \sqrt{ b^2+(a^{-1}+a)^2}} & \frac{b}{\pm \sqrt{ b^2+(a^{-1}+a)^2}} \\
-\frac{b}{\pm \sqrt{ b^2+(a^{-1}+a)^2}} & \frac{a^{-1}+a}{\pm \sqrt{ b^2+(a^{-1}+a)^2}} \end{array} \right), \]
 choosing $\hbox{sign}(\pm \sqrt{ b^2+(a^{-1}+a)^2})=\hbox{sign}\  b$ if 
$b \neq 0$ and $+1$ for $b=0$. 
The section $\sigma _1$ corresponds to the hyperbolic plane loop 
(cf. \cite{loops}, pp. 283-284).
\newline
Since $\big [ [\bf{m},\bf{m}], \bf{m} \big ] \subseteq \bf{m}$ and each element $g \in G$ can 
uniquely be represented as a product  
$g= m h$ with $m \in \exp {\bf m}$ and $h \in H_3$ we have a global Bol loop $L$ 
defined on the factor space $G/H_3$ (cf. 
\cite{kloops},  Corollary 3.11, p. 51 and \cite{loops}, Lemma 1.3, p. 17 ).
In $L$ there is a normal subgroup $\tilde{G}$ isomorphic to $SO_2(\mathbb R)$ 
and the factor loop $L/\tilde{G}$  is isomorphic to the hyperbolic plane 
loop. Therefore is $L$ the unique Scheerer extension of the Lie group $SO_2(\mathbb R)$ by the 
hyperbolic plane loop (cf. \cite{loops}, Section 2). This loop is a left A-loop, 
because of $[ {\bf h}, {\bf m}] \subseteq {\bf m}$. But it is not a Bruck loop since there is no involutory automorphism $\sigma :{\bf g} \to {\bf g}$ such that $\sigma ({\bf m})=- {\bf m}$ and $\sigma ({\bf h})={\bf h}$.   

\begin{Theo} There are precisely three isotopism classes ${\cal C}_1$, 
${\cal C}_2$, ${\cal C}_3$ of  connected differentiable Bol loops $L$ with dimension $3$ 
such that the group $G$ generated by the left translations $\{ \lambda_x ; x 
\in L\}$ is a $4$-dimensional non-solvable Lie group.
Every loop in the class ${\cal C}_1$, respectively ${\cal C}_2$ is 
the direct 
product of a $2$-dimensional loop isotopic to the hyperbolic plane loop 
 with the Lie group 
$\mathbb R $, respectively $SO_2(\mathbb R)$. In the class ${\cal C}_3$ is contained up to 
isomorphisms only the Scheerer extension of the Lie group $S0_2(\mathbb R)$ by the 
hyperbolic plane loop. 
\end{Theo}

\section{3-dimensional Bol loops belonging to  
\newline 5-dimensional non-solvable Lie groups}
In this section we seek for $3$-dimensional connected differentiable global Bol loops having a $5$-dimensional non-solvable Lie group $G$ as the group 
topologically generated by
their left translations. 
\newline
Any $5$-dimensional non-solvable Lie group is locally isomorphic to one of 
the following groups:
\newline
1)\ $PSL_2(\mathbb R) \times \mathbb R^2$
\newline
2)\ $PSL_2(\mathbb R) \times$ $\cal{L}$$_2$, where $\cal{L}$$_2 \cong 
\{x \mapsto ax +b; a>0, b \in \mathbb R \}$.     
\newline
3)\ $SO_3(\mathbb R) \times \mathbb R^2$
\newline
4)\ $SO_3(\mathbb R) \times$ $\cal {L}$$_2$
\newline
5) The semidirect product $G=SL_2(\mathbb R) \ltimes \mathbb R^2$ such that $SL_2(\mathbb R)$ 
acts in the natural way on 
$\mathbb R^2$. The group $G$ may be seen as the connected component of the group of area 
preserving affinities of the real plane. 
\newline
Since there is no epimorphism 
from a $6$-dimensional Lie algebra ${\bf g}^*$ given in section 3 onto one of the $5$-dimensional Lie 
algebras ${\bf g}$  listed above only the Lie algebra 
${\bf g}=sl_2(\mathbb R) \ltimes \mathbb R^2$, which is the Lie algebra of the isometry group 
of a symmetric space (see case {\bf 7} in section 3) can occur as the group topologically 
generated by the left translations of a 
$3$-dimensional connected differentiable Bol loop.  
Every $2$-dimensional connected subgroup $H$ of $G$, which does not have any non-trivial 
normal subgroup 
of $G$, contains $1$-parameter subgroups $S$ fixing a point or contains a $1$-dimensional 
translation group. The Lie algebra of $S$ has as generator an elliptic, parabolic, or 
hyperbolic element of ${\bf g}$.  
Using the relations of {\bf 7} we see that ${\bf m}$ contains  elements  
corresponding to translations and ${\bf m} \cap sl_2(\mathbb R)$ has with $e_2$ a hyperbolic, 
with $e_3$ an elliptic and with $e_2+e_3$ a parabolic element.  
Therefore $\exp {\bf m}$ contains elements conjugate to elements of $H$. 
According to Lemma 3 there is  no $3$-dimensional 
Bol loop $L$ as section in the group $PSL_2(\mathbb R) \ltimes \mathbb R^2$. 
This consideration yields the following
\begin{Theo} There is no $3$-dimensional connected differentiable 
global Bol loop $L$ having a $5$-dimensional non-solvable Lie group as the 
group topologically generated by its left translations.
\end{Theo}

\section{3-dimensional Bol loops with 6-dimensional non-solvable Lie 
groups as their left translations groups}

A $6$-dimensional non-solvable Lie group $G$ can be occur as the group topologically generated by the left translations of a $3$-dimensional connected differentiable Bol loop, if $G$ is isomorphic to a $6$-dimensional Lie group $G^*$ listed in the section 3. Then we have to deal with the cases, where $G$ is locally isomorphic to one of the following Lie groups: 
\newline
1) The connected component of the euclidean motion group on $\mathbb R^3$. 
\newline
2) The semidirect product of $PSL_2(\mathbb R) \ltimes \mathbb R^3$, where the action of 
$PSL_2(\mathbb R)$ on $\mathbb R^3$ is the adjoint action of $PSL_2(\mathbb R)$ on its Lie algebra. 
\newline
Let us consider the first case. 
The $3$-dimensional subgroups $H$ of $G=SO_3(\mathbb R) \ltimes \mathbb R^3$,  
which does not contain any non-trivial normal subgroup of $G$, 
are locally isomorphic to the following 
subgroups: 
Either $H$ has the shape as in case c)
of Proposition 2 and there is no $3$-dimensional differentiable Bol loop corresponding to the pair 
$(G,H)$ 
or $H=\{(a,0), a \in SO_3(\mathbb R) \}$. 
We consider the $\mathbb R$-basis {\bf B3} of the Lie algebra ${\bf g}=so_3(\mathbb R) \ltimes \mathbb R^3$ given in the case {\bf 6.1} (Section 3). 
Then the Lie algebra $\bf{h}$ of $H$ is generated by the 
 basis elements $e_2, e_3, e_4$ and the Lie triple system ${\bf m}$ 
is generated by the basis elements $e_1, e_2, e_3$.  
Then ${\bf h} \cap {\bf m}=\langle e_2, e_3 \rangle$ which is  
a contradiction.  
\newline
\newline
Now we discuss the second case. Then the group multiplication is 
given by 
\[ (A_1, X_1) \circ (A_2, X_2)= (A_1\ A_2,\ A_2^{-1}\ X_1\ A_2+ X_2), \]
where $(A_i,X_i)$, $i=1,2$ are  two elements of $G$ such that $X_i$ $i=1,2$ 
are represented by $2 \times 2$ real  matrices with trace $0$. We use the 
$\mathbb R$-basis {\bf B4} of the Lie algebra ${\bf g}=sl_2(\mathbb R) \ltimes \mathbb R^3$ ({\bf 6.3}, Section 3). The $3$-dimensional subalgebras ${\bf h}$ of ${\bf g}$ containing no non-zero ideal of ${\bf g}$ are the following: 
\newline
a) $\langle e_2,\ e_5,\ e_1+e_6 \rangle $,
\newline
b) $ \langle e_2+k\ e_5,\ e_1,\ e_6 \rangle$, where $k \in \mathbb R$,
\newline
c) $\langle e_3+e_4,\ e_5,\ e_1-e_6 \rangle $,
\newline 
d) $\langle e_2,\ e_3+e_4,\ e_1-e_6 \rangle $,    
\newline
e) $\langle e_2,e_3,e_4 \rangle $,
\newline
f) $\langle e_4,e_5,e_6 \rangle $.     
\newline
First we deal with the Lie triple system ${\bf m}$ of ({\bf 6.2} in section 3) generated by the 
basis elements 
$e_4, e_2, e_6$. If the Lie algebra ${\bf h}$ has the shape a), b), d), e) and f) then the intersection of  ${\bf h}$
with ${\bf m}$ are the subspaces 
$\langle e_2 \rangle$, 
$\langle e_6 \rangle$, $\langle e_2 \rangle$, $\langle e_2,e_4 \rangle$, or 
$\langle e_4,e_6 \rangle$ respectively.  
In the case c)  the elements $e_2+e_4 \in {\bf m}$ and $e_3+e_4 \in {\bf h}$
are parabolic and the parabolic elements of $sl_2(\mathbb R)$ are conjugate to 
each other. This is a contradiction to Lemma 3.

\bigskip
\noindent
Let us treat the Lie triple system ${\bf m}$ of ({\bf 6.3} in section 3) which is generated by 
the basis elements $e_1, e_2, e_3$ and the image of ${\bf m}$ under the exponential map. 
\newline
The exponential map $\exp : {\bf g} \to G$ is defined in the following way: For
$X \in {\bf g}$ we have $\exp \ X= \gamma _X(1)$, where $\gamma _X(t)$ is the 
1-parameter subgroup of $G$ with the property 
$ \frac{d}{dt} \big |_{t=0} \gamma _X(t)=X$.
In the 1-parameter subgroup $ \alpha (t)=( \beta (t), \gamma (t))$ of $G$
with the conditions $\alpha (t=0)=(1,0)$ and   
$ \frac{d}{dt} \big |_{t=0} \alpha (t)=(X_1,X_2)=X \in \bf{g}$ the first component 
$\beta (t)$ is the 1-parameter subgroup of 
$PSL_2( \mathbb R)$, and the second component satisfies 
\newline
\centerline{$\displaystyle \frac{d}{dt} \gamma (t)= \frac{d}{ds} \big| _{s=0} \gamma (t+s)=
- \frac{d}{ds} \big| _{s=0} \beta (s) \gamma(t)+\gamma(t)  
\frac{d}{ds} \big| _{s=0} \beta (s)+ \frac{d}{ds} \big| _{s=0} \gamma (s)=$} 
\[ -X_1 \gamma(t) + \gamma(t) X_1 +X_2. \] 
For $X_1= \left( \begin{array}{rr}
a & b \\
c & -a \end{array} \right)$, and $X_2= \left( \begin{array}{rr}
k & u \\
y & -k \end{array} \right)$, with $a,b,c,k,u,y \in \mathbb R$ 
 and 
 $\gamma (t)=\left( \begin{array}{rr}
r(t) & s(t) \\
v(t) & -r(t) \end{array} \right)$,
one has \[ \frac{d}{dt} \gamma (t)= \left( \begin{array}{rr}
 \frac{d}{dt} r(t)  & \frac{d}{dt} s(t) \\
 \frac{d}{dt} v(t)  & -\frac{d}{dt} r(t)  \end{array} \right)=
\left( \begin{array}{rr}
-a & -b \\
-c & a \end{array} \right) \left( \begin{array}{rr}
r(t) & s(t) \\
v(t) & -r(t) \end{array} \right)+ \] \[ \left( \begin{array}{rr}
r(t) & s(t) \\
v(t) & -r(t) \end{array} \right) \left( \begin{array}{rr}
a & b \\
c & -a \end{array} \right)+ \left( \begin{array}{rr}
k & u \\
y & -k \end{array} \right) \]
 with the following properties: 
\[ r(0)=s(0)=v(0)=0,\  \frac{d}{dt} \big| _{t=0} r(t)=k, \ 
\frac{d}{dt} \big| _{t=0} s(t)=u,\ \frac{d}{dt} \big| _{t=0} v(t)=y,\] 
\[ \frac{d}{dt} r(t)=-b v(t)+c s(t)+k,\ \frac{d}{dt} s(t)=2(b r(t)+a s(t))
+u,\] 
\[ \frac{d}{dt} v(t)=2(a v(t)-c r(t))+y.\] The solution of this inhomogeneous 
system of linear differential equations is:  
\newline
\[r(t)= \frac{1}{8(a^2+bc)^{\frac{3}{2}}} (e^{2 
\sqrt {a^2+bc}\  t}- e^{-2 \sqrt {a^2+bc}\  t} )(-acu-bay+2kcb) +\] \[ 
\frac{1}{8(a^2+bc)} \big [ (e^{ 
\sqrt {a^2+bc}\  t}- e^{- \sqrt {a^2+bc}\  t})^2 (-cu+by) + 
\  t (8ka^2+4acu+4aby) \big], \] 
\[s(t)= \frac{1}{8(a^2+bc)^{\frac{3}{2}}} (e^{2 
\sqrt {a^2+bc}\  t}- e^{-2 \sqrt {a^2+bc}\  t} )(-b^2y+ubc-2bak+2a^2u) +\] 
\[ \frac{1}{8(a^2+bc)} \big[ (e^{ \sqrt {a^2+bc}\  t}- e^{- \sqrt {a^2+bc}
\  t})^2(-2bk +2au) +  t (4b^2y+4ubc+8kab) \big], \]
\[v(t)= \frac{1}{8(a^2+bc)^{\frac{3}{2}}}  (e^{2 
\sqrt {a^2+bc}\  t}- e^{-2 \sqrt {a^2+bc}\  t} )(2ya^2-2cka+bcy-c^2u) +\] \[
\frac{1}{8(a^2+bc)} \big[ (e^{ \sqrt {a^2+bc}\  t}- e^{- \sqrt {a^2+bc}\  t})^2  (-2ay+2ck) + t (8cka+4c^2u+4bcy) \big] . \]
It follows  $\exp {\bf m}=\{ \exp (X_1,X_2) \}$, where  
$X_1 \in \left \{ \left( \begin{array}{cc}
\lambda _2 & \lambda_3  \\
 \lambda _3 & -\lambda _2 \end{array} \right), \lambda_2, \lambda_3 \in \mathbb R \right \}$ and 
$X_2 \in \left  \{ \left( \begin{array}{cc}
0 & -\lambda _1 \\
\lambda _1 & 0
\end{array} \right), \lambda_1 \in \mathbb R \right \}$.
\newline
The first component of $ \exp {\bf m}$ is 
\[ (\exp (X_1,X_2))_1= \left( \begin{array}{cc}
\cosh \sqrt {A}+ \displaystyle \frac{ \sinh \sqrt {A}}{\sqrt {A}} \lambda_2 & 
\displaystyle \frac{ \sinh \sqrt {A}}{\sqrt {A}} \lambda _3 \\
\displaystyle \frac{ \sinh \sqrt {A}}{\sqrt {A}}\lambda _3  
& \cosh \sqrt {A}-\displaystyle \frac{ \sinh \sqrt {A}}{\sqrt {A}} 
\lambda_2 \end{array} 
\right), \] where 
$A=\lambda _2^2+\lambda _3^2$ (cf. \cite{old}, p. 24);
the second component of $\exp {\bf m}$ is 
$(\exp (X_1,X_2))_2=\left( \begin{array}{rr}
r(1) & s(1) \\
v(1) & -r(1) \end{array} \right) $, where 
\[ r(1)=\frac{\lambda _3\  \lambda _1}{4(\lambda _2^2 + \lambda _3^2)}\ 
(e^{ \sqrt {\lambda _2^2 + \lambda _3^2}}- e^{ - \sqrt 
{\lambda _2^2 + \lambda _3^2}})^2, \] \[ s(1)= \frac{- \lambda _1}
{4 \sqrt {\lambda _2^2 + \lambda _3^2} }(e^{2 \sqrt {\lambda _2^2 + 
\lambda _3^2}}- e^{ -2 \sqrt 
{\lambda _2^2 + \lambda _3^2}}) -\frac{\lambda _2 \lambda _1}
{4(\lambda _2^2 + \lambda _3^2)}
(e^{ \sqrt {\lambda _2^2 + \lambda _3^2}}- e^{ - \sqrt 
{\lambda _2^2 + \lambda _3^2}})^2, \]
\[ v(1)= \frac{\lambda _1}
{4 \sqrt {\lambda _2^2 + \lambda _3^2} }(e^{2 \sqrt {\lambda _2^2 + 
\lambda _3^2}}- e^{ -2 \sqrt 
{\lambda _2^2 + \lambda _3^2}}) -\frac{\lambda _2 \lambda _1}
{4(\lambda _2^2 + \lambda _3^2)} 
(e^{ \sqrt {\lambda _2^2 + \lambda _3^2}}- e^{ - \sqrt 
{\lambda _2^2 + \lambda _3^2}})^2. \] 
The Lie algebra ${\bf h}$ of the stabilizer cannot have the shape a), b), d) and e) 
since in these cases the intersection of  ${\bf h}$
with ${\bf m}$ is the subspace $\langle e_2 \rangle$,  $\langle e_1 \rangle$,
$\langle e_2 \rangle$ and  $\langle e_2,e_3 \rangle$ respectively. 
\newline
In the case c) we have ${\bf h} \cap {\bf m}= 0 $, 
there is no element of ${\bf m}$ which is conjugate to an  element of 
${\bf h}$ and 
\newline
\centerline{
$H=\left \{ \left( \left( \begin{array}{cc}
1 & b \\
0 & 1 \end{array} \right) , \left( \begin{array}{cc}
e & f \\
0 & -e \end{array} \right) \right); b,e,f \in \mathbb R  \right 
\}$.}
To prove that in this case there is no global 
Bol loop we may proceed as we did it for the group $G=PSL_2(\mathbb R) \times \mathbb R$ and 
for the stabilizer 
\newline
$\exp {\bf h}_2 = 
\left \{ \left ( \left ( \begin{array}{cc} 
1 & b \\
0 & 1 \end{array} \right ), b \right); b \in \mathbb R \right \}$ since the first component of $H$ 
is equal to the first component of  $\exp {\bf h}_2$ and the first components of $\exp {\bf m}$ 
are the same in both cases. 
\newline
\newline
Now we assume that the stabilizer $H$ has the shape f).  
Since the group $SL_2(\mathbb R)$ has no $3$-dimensional linear representation the group $G$ is 
isomorphic to the 
semidirect 
product of $PSL_2( \mathbb R) \ltimes \mathbb R^3$ and $H$ may be chosen as   
\newline
$H=\left \{ \left( \pm \left( \begin{array}{rr}
\cos {t} & \sin {t} \\
- \sin {t} & \cos {t} \end{array} \right),  \left( \begin{array}{rr}
-x & y \\
y  & x  \end{array} \right) \right), t \in [0, 2 \pi ), x,y \in \mathbb R 
\right \}$. 
\newline
We have $[{\bf m}, {\bf m}]= {\bf h}$, ${\bf g}={\bf m} \oplus [ {\bf m}, {\bf m}]$ and 
therefore the corresponding loop 
$L$ realized on $G/H$ is a Bruck loop. We prove that this loop is global.  
\newline
For the connected simple 
Lie group 
$PSL_2( \mathbb R)$ there exists a unique decomposition 
$g=\left( \begin{array}{cc}
a & 0 \\
b & a^{-1} \end{array} \right) \left(\pm \left( \begin{array}{rr}
\cos {t} & \sin {t} \\
- \sin {t} & \cos {t} \end{array} \right) \right)$, where $a>0, b \in \mathbb R$, 
\newline
$t \in [0, 2 \pi )$ (cf. \cite{Iwasawa}, p. 525). 
Since 
\[ \left( \pm \left( \begin{array}{rr}
\cos {t} & -\sin {t} \\
 \sin {t} & \cos {t} \end{array} \right) \right) \left( \begin{array}{rr}
0 & u \\
-u & 0 \end{array} \right) \left( \pm \left( \begin{array}{rr}
\cos {t} & \sin {t} \\
- \sin {t} & \cos {t} \end{array} \right) \right)+ \left( \begin{array}{rr}
k & l \\
l & -k \end{array} \right)=\] \[  \left( \begin{array}{rr}
0 & u \\
-u & 0 \end{array} \right)+  \left( \begin{array}{rr}
k & l \\
l & -k \end{array} \right)=  \left( \begin{array}{cc}
k & u+l \\
l-u & -k \end{array} \right), \] each element of $G$ can uniquely be written as
$\left( \pm \left( \begin{array}{cc}
a & b \\
c & d  \end{array} \right), \left( \begin{array}{rr}
x & y \\
z & -x \end{array} \right) \right)= \left( \left( \begin{array}{ll}
a_1 & 0 \\
b_1 & a_1^{-1} \end{array} \right),  \left( \begin{array}{rr}
0 & u \\
-u & 0 \end{array} \right) \right) \circ  \left( \pm \left( \begin{array}{rr}
\cos {t} & \sin {t} \\
- \sin {t} & \cos {t} \end{array} \right), \left( \begin{array}{rr}
k & l \\
l & -k \end{array} \right) \right)$ with 
\newline 
$a,b,c,d \in \mathbb R ,\ a d-b c=1,\ x,y,z \in \mathbb R, a_1>0, b_1,u,k,l \in \mathbb R,\ t \in [0,2 \pi)$,
such that $k=x, l=\displaystyle \frac{y+z}{2}, u=\displaystyle \frac{y-z}{2}$. 
The submanifold $\exp {\bf m}$ is the image of a section 
$\tau: G/H \to G$ if and only if each element $g \in G$ can uniquely be 
represented as a product $g=m\ h$ with $m \in \exp {\bf m}$ and 
$h \in H$ or equivalently $m=g\ h^{-1}$. It is sufficient to prove this for each element $g \in G$ with the shape 
\[ \left \{ \left(  \left( \begin{array}{ll}
a & 0 \\
b & a^{-1} \end{array} \right),  \left( \begin{array}{rr}
0 & u \\
-u & 0 \end{array} \right) \right), a>0,b,u \in \mathbb R \right \}. \] 
The first component of $\exp {\bf m}$ is precisely  
the section $\sigma _1$ of the hyperbolic plane loop given before Theorem 6. Therefore  
for given $a>0,b \in \mathbb R$ we have unique 
$ \lambda_2, \lambda_3 \in \mathbb R, t \in [0, 2 \pi)$ such that 
\[ \left( \begin{array}{cc}
 \cosh \sqrt {A}+  \frac{ \sinh \sqrt {A}}{\sqrt {A}} \lambda_2 & 
 \frac{ \sinh \sqrt {A}}{\sqrt {A}} \lambda _3 \\
\frac{ \sinh \sqrt {A}}{\sqrt {A}}\lambda _3  
& \cosh \sqrt {A}- \frac{ \sinh \sqrt {A}}{\sqrt {A}} 
\lambda_2 \end{array} 
\right)= \] 
\[ \left( \begin{array}{ll}
a & 0 \\
b & a^{-1} \end{array} \right) \left( \begin{array}{rr}
\cos {t} & \sin {t} \\
- \sin {t} & \cos {t} \end{array} \right), \] 
where $A=\lambda _2^2+\lambda _3^2$. Hence we have to consider the second component 
 of $\exp {\bf m}$.  
For given $u, \lambda_2, \lambda_3$ we 
have to find unique $\lambda_1,k,l \in \mathbb R$ such that 
\[ \left( \begin{array}{rr}
r(1) & s(1) \\
v(1) & -r(1) \end{array} \right)= \left( \begin{array}{cc}
k & l+u \\
l-u & -k \end{array} \right), \]
where 
$r(1), s(1), v(1)$ are values of functions, which depend on 
the variables 
$\lambda_1, \lambda_2,\lambda_3$. 
Since for $\lambda_1$ we obtain the equation 
\[ 2u= \displaystyle \frac{-\lambda_1}{2 \sqrt{\lambda_2^2+\lambda_3^2}} 
(e^{2 \sqrt{ \lambda_2^2+\lambda_3^2}}-e^{-2 \sqrt{\lambda_2^2+\lambda_3^2}}) 
\] we have the unique solutions 
$\lambda_1=\displaystyle \frac{-4 u \sqrt {\lambda_2^2+\lambda_3^2}}
{\displaystyle 
e^{2 \sqrt{ \lambda_2^2+\lambda_3^2}}-e^{-2 \sqrt{\lambda_2^2+\lambda_3^2}}}$,   $k=r(1)$ and $l=\displaystyle \frac {s(1)+v(1)}{2}$.
\newline
Summarizing the above discussion we see that $\exp {\bf m}$ is an image of a global section 
$\sigma:G/H \to G$. The section $\sigma $ determines a global Bol loop $({\hat L}, \ast )$ since 
$\big [ [\bf{m},\bf{m}],\bf{m} \big ] \subseteq \bf{m}$ and therefore the equation $x \ast a=b$ has the 
solution $x= a^{-1} \ast [ (a \ast b) \ast a^{-1}]$ for all $a,b \in {\hat L}$ 
(cf. \cite{kloops},  Corollary 3.11, p. 51 and \cite{loops}, Lemma 1.3, p. 17 ). 
We call the loop $({\hat L}, \ast )$ the pseudo-euclidean space loop. 
\newline
\newline
In order to determine within the isotopism class ${\cal C}$ of the pseudo-euclidean space loop 
$({\hat L}, \ast )$ the isomorphism classes we consider arbitrary complements ${\bf m}$ to the Lie algebra ${\bf h}$ 
of the stabilizer $H$ of $e \in L$ in ${\bf g}$. These complements have the form: 
\[ \langle e_1+a_1 e_4+a_2 e_5+a_3 e_6, e_2+b_1 e_4+b_2 e_5+b_3 e_6, 
e_3+c_1 e_4+c_2 e_5+c_3 e_6 \rangle \] where  
$a_1,a_2,a_3,b_1,b_2,b_3,c_1,c_2,c_3 \in \mathbb R$. 
Such complements are Bol complements if and only if  
\[ {\bf m}={\bf m}_{b_3,c_3,c_2}= \langle e_1-c_3 e_6+b_3 e_5, e_2+c_2 e_6+b_3 e_4, e_3+c_2 e_5+c_3 e_4 \rangle ,\]
 where 
$c_2,c_3,b_3 \in \mathbb R$, and $b_3^2+c_3^2 \neq 1$. 
According to remark in section 2 the subspace ${\bf m}_{b_3,c_3,c_2}$ belongs to a loop within 
the isotopism class ${\cal C}$ 
if and only if the subalgebra ${\bf h}_{b_3,c_3,c_2}= [{\bf m}_{b_3,c_3,c_2}, {\bf m}_{b_3,c_3,c_2}]$   
is the Lie algebra of the semidirect product of a $2$-dimensional normal translation group $T$ 
by a $1$-dimensional rotation group $S$. It follows that then ${\bf h}_{b_3,c_3,c_2}={\bf t}  
\oplus {\bf s}$, where ${\bf t}$ respectively ${\bf s}$ is the Lie algebra of $T$ respectively 
of $S$ and 
the Cartan-Killing form $k$ given by {\bf K2} in section 3 is zero on ${\bf t}$ and negative 
on ${\bf s}$. 
Since $[{\bf m}_{b_3,c_3,c_2}, {\bf m}_{b_3,c_3,c_2}]=\langle d_1, d_2, d_3 \rangle$ with 
\newline 
\centerline{$d_1=(1-b_3^2)e_6-c_3 e_1 -c_3 b_3 e_5$, $d_2=(1-c_3^2)e_5 +b_3 e_1-c_3 b_3 e_6$,} 
\centerline{$d_3=e_4+c_3 e_3+b_3 e_2 +c_2 c_3 e_5 +b_3 c_2 e_6$,} the subalgebra 
$\langle d_1, d_2 \rangle$ is an ideal  
 of ${\bf h}_{b_3,c_3,c_2}$ and $k(d_3) < 0$ if and only if $b_3^2+c_3^2<1$. 
\newline
The Bol complements ${\bf m}_{b_3,c_3,c_2}$ with $b_3^2+c_3^2>1$ are according to Remark in section 2 
symmetric spaces  belonging to the Lie triple system given 
in {\bf 6.2} (section 3). But the local Bol loop corresponding to 
${\bf g}={\bf m}_{b_3,c_3,c_2} \oplus {\bf h}$ cannot be extended to a global Bol loop.        
\newline
The loop $L$  belonging to $(G,H, \exp {\bf m}_{b_3,c_3,c_2})$, $b_3^2+c_3^2<1$ and the loop 
$L'$
corresponding to $(G,H, \exp {\bf m}_{b_3',c_3',c_2'})$, $b_3'^2+c_3'^2<1$   
are isomorphic if and only if there exists an automorphism $\alpha $ of 
${\bf g}$ such that $\alpha ({\bf m}_{b_3',c_3',c_2'})=
{\bf m}_{b_3,c_3,c_2}$ 
and 
$\alpha ({\bf h})={\bf h}$. The automorphism   
group $A$ of the Lie algebra ${\bf g}$ leaving ${\bf h}$ invariant consists 
of 
linear mappings determined by the basis transformations:  
\newline
 $ \alpha \  (e_1)=a\  e_1$, 
\newline
$\alpha \  (e_6)=b_2\ e_6 + b_4\  e_5$,
\newline
$\alpha \  (e_5)=- \varepsilon\ b_4\  e_6\ + \varepsilon \  b_2\ e_5$,
\newline
$\alpha \ (e_4)=\varepsilon \ e_4 + d_5\ e_6 + d_6\ e_5$,
\newline
$\alpha \ (e_2)=\displaystyle \frac{b_2}{a}\ e_2 
+ \frac{b_4}{a}\ e_3 + \left ( \frac{ - \varepsilon \  
d_5\ b_4\ + \varepsilon \  d_6\ b_2}{a} \right)\ e_1
+ f_5\ e_6 +f_6\ e_5 $,
\newline
$\alpha \ (e_3)=\displaystyle \varepsilon \frac{b_2}{a}\ e_3
- \varepsilon \ \frac{b_4}{a}\ e_2+ \left (
\frac{-\ d_5\ b_2\ -\ d_6\ b_4}{a} \right ) e_1 - \varepsilon f_6\ e_6\ + 
\varepsilon \ f_5\ e_5 $, 
\newline
where $\varepsilon \in \{1,-1 \}$ and for the real parameters 
$a,b_2,b_4,d_5,d_6,f_5,f_6$ the properties $a^2=b_2^2+b_4^2 \neq 0$ and $b_2\ f_6= b_4\ f_5$ are 
satisfied. 
\newline 
According to {\bf K2} (Section 3) the Cartan-Killing 
form $k$ is invariant under $\alpha $. Hence the $3$ angles 
between ${\bf m}_{b_3,c_3,c_2}$ and $\hat{{\bf m}}:={\bf m}_{0,0,0}$ respectively 
between 
${\bf m}_{b_3',c_3',c_2'}$ and 
$\hat{{\bf m}}$ are equal with respect to $k$. 
\newline
Now we compute the isomorphism classes of the loops 
$L_{b_3,c_3,c_2}$.  Since the automorphism $\gamma $ of ${\bf g}$ given by 
$\gamma (e_1)=-e_1$, $\gamma (e_2)=-e_2 -2 c_2 \ e_4$, $\gamma (e_3)=- e_3 -2 c_2 \ e_5$, 
$\gamma (e_i)=e_i$ for $i=4,5,6$ maps ${\bf m}_{b_3,c_3,c_2}$ onto ${\bf m}_{b_3,c_3,0}$   
we seek for which real triples $(b_3',c_3',0)$ and $(b_3,c_3,0)$ there is an automorphism 
$\alpha $ in $A$, such that 
$\alpha ( {\bf m}_{b_3',c_3',0} )={\bf m}_{b_3,c_3,0}$.
This condition is equivalent to the following 
 system of equations: 
\newline
$\begin{array}{rcl}
-a\ c_3 &= & -c_3'\ b_2 - \varepsilon \ b_3'\ b_4,\\ 
a\ b_3 &= & -c_3'\ b_4 + \varepsilon \ b_3'\ b_2, \\  
\displaystyle a\ (f_5\ +b_3'\ d_5) & = &  
(\varepsilon\  d_5\ b_4\ -\varepsilon\ d_6\ b_2)\ c_3, \\
a\ (f_6+b_3'\ d_6) & = &  
(- \varepsilon\ d_5\ b_4\ + \varepsilon \ 
d_6\ b_2)\ b_3, \\
a\ (- \varepsilon\ f_6 +c_3'\ d_5) & = & \displaystyle 
(d_5\ b_2\ +d_6\ b_4)\ c_3, \\
\displaystyle a\ (\varepsilon\ f_5+c_3'\ d_6) & = &
(-d_5\ b_2\ -d_6\ b_4)\ b_3. \end{array}$
\newline
Solving this system of equations we see that the loop 
$\hat{L_{d}}=L_{d,0,0}$ for $0 \le d <1$ is a representative 
of the isomorphism class consisting of loops $L_{b_3,c_3,c_2}$ with $\sqrt {b_3^2+c_3^2}=d$. 
The pseudo-euclidean space loop $ÿ({\hat L}, \ast )=L_{0,0,0}$ is a representative of the isomorphism class 
consisting of the loops $L_{0,0,c_2}$, $c_2 \in \mathbb R$. 
\newline
A loop isomorphic to ${\hat L}$ is a Bruck loop and hence ${\bf m}_{0,0,c_2}$ 
is reductive to ${\bf h}$ in ${\bf g}$. If ${\bf m}_{d,0,0}$ with $d \neq 0$ were reductive to 
${\bf h}$ in the Lie algebra ${\bf g}$ then for $e_4 \in {\bf h}$ and $e_1+d\ e_5 \in {\bf m}_{d,0,0}$ 
the element $[e_4, e_1+d\ e_5 ]=d\ e_6 \in
{\bf h}$ would be contained in ${\bf m}_{d,0,0}$, which is a contradiction.  
\newline
\newline
By $E(2,1)$ we denote the pseudo-euclidean space the points of which are represented by 
the matrices 
\[(I , Y)=\left( I , \left( \begin{array}{cc}
x & k+l \\
k-l & -x \end{array} \right) \right), \] where $x,k,l \in \mathbb R$, and which has as norm 
$|| (I, Y) ||=x^2+k^2-l^2$.  
The group $G$ acts on the space $E(2,1)$ in the following way:
For given $(A, X) \in G$ and $(I, Y) \in E(2,1)$
\[\ \ (\ast )\quad \quad (A,X) \ast (I,Y) = (I ,\ A^{-1} Y A+ X). \]
The norm is invariant under the action of $G$, therefore  $G$ is the 
connected component of the motion group of $E(2,1)$.  
\newline
The $3$-dimensional pseudo-euclidean geometry $E(2,1)$ has also a representation 
${\cal R}$ on the affine space $\mathbb R^3$ such that the motion group 
consists of the affine mappings 
\[ \ \ (B,b): (x,y,z) \mapsto (x,y,z) B^T + b,  \]
where $B=\left ( \begin{array}{ccc} 
a_1 & a_2 & a_3 \\
b_1 & b_2 & b_3 \\
c_1 & c_2 & c_3 \end{array} \right )$ with $\det B=1$, $a_1^2+a_2^2=a_3^2+1$,  
$b_1^2+b_2^2=b_3^2+1$, $c_1^2+c_2^2=c_3^2-1$, and $b=(b_1, b_2, b_3)$ 
(\cite{benz}, Kapitel 6).   
The mappings 
\newline
\centerline{ $\omega : \left( I , \left( \begin{array}{cc}
k & l+n \\
l-n & -k \end{array} \right) \right) \mapsto (k,l,n)$} 
and 
\centerline{ $\Omega : \left ( \pm \left ( \begin{array}{cc}
a & b \\
c & d \end{array} \right ), \left ( \begin{array}{cc}
x & y+z \\
y-z & -x \end{array} \right ) \right) \mapsto $}
\centerline{ $\left ( \left ( \begin{array}{ccc} 
da+bc & cd-ba & cd+ba \\
bd-ca & \displaystyle \frac{a^2+d^2-b^2-c^2}{2} & \displaystyle \frac{d^2+b^2-a^2-c^2}{2} \\
bd+ca & \displaystyle \frac{d^2-b^2-a^2+c^2}{2} & \displaystyle \frac{d^2+b^2+a^2+c^2}{2} 
\end{array} \right ), ( x, y, z ) \right)$} 
establish an isometry from $E(2,1)$ onto ${\cal R}$ (\cite{klein}, pp. 97-103). 
\newline
The stabilizer 
$H$ which is the image of the Lie algebra of the shape f) under the exponential map leaves 
in $E(2,1)$ the plane $P$ 
consisting of the points 
$\left \{  \left(I, \left( 
\begin{array}{rr}
x & y \\
y & -x \end{array} \right) \right), x,y \in \mathbb R \right \}$ invariant. The points of $P$ 
satisfy  $x^2+y^2>0$. The planes of 
${\cal R}$ the points $(x,y,z)$ of which satisfy $x^2+y^2>z^2$ are called euclidean planes. 
The connected component $\Omega (G)$ of the motion group of  ${\cal R}$  acts transitively on 
the set $\Psi $ of the euclidean planes and the set $\Omega (\exp {\bf m}_{0,0,0})$ is sharply 
transitive on $\Psi $. The planes of $\Psi $ can be taken as the points of the pseudo-euclidean 
space loop $({\hat L}, \ast )$ such that the multiplication is given by 
\newline
\centerline{$(\ast \ast) \quad \quad Q_1 \ast Q_2= \tau _{P,Q_1} (Q_2)$, \quad for all 
$Q_1, Q_2 \in \Psi $,  } 
where $\tau_{P,Q_1}$ is the unique element of $\Omega (\exp {\bf m}_{0,0,0})$ mapping 
the plane $P$, which is the identity of ${\hat L}$ onto $Q_1$.

\bigskip
\noindent
From the above discussion we obtain the 
\begin{Theo} There is only one isotopism class $\cal C$ of the $3$-dimensional
connected differentiable Bol loops $L$ such that the group $G$ generated by 
the left translations $\{ \lambda_x; x \in L \}$ is a $6$-dimensional 
non-semisimple and non-solvable Lie group.
The group $G$ is isomorphic to $PSL_2(\mathbb R) \ltimes \mathbb R^3$, where 
the action of $PSL_2(\mathbb R)$ on $\mathbb R^3$ is the adjoint action of 
$PSL_2(\mathbb R)$ on its Lie algebra. 
The stabilizer $H$ of 
$e \in L$ in $G$ is isomorphic to  
\[H=\left \{ \left(\pm \left( \begin{array}{rr}
\cos t & \sin t \\
-\sin t & \cos t \end{array} \right) , \left( \begin{array}{rr}
-x & y \\
y & x \end{array} \right) \right); t \in [0 , 2 \pi ],x,y \in \mathbb R 
\right \}. \] 
The class ${\cal C}$ consists of Bol loops  $L_{b_3,c_3,c_2}$, where 
$b_3,c_3,c_2$ are real parameters and $b_3^2+c_3^2<1$.
As a representative of ${\cal C}$ one may choose the pseudo-euclidean space loop 
$\hat{L}=L_{0,0,0}$. This loop is a Bruck loop. 
\newline
The tangential spaces $T_1 \sigma (G/H)$ of the loops $L_{b_3,c_3,c_2}$ at 
$1 \in G$  
are the Bol complements ${\bf m}_{b_3,c_3,c_2}$.
The loops $L$ and $L'$ in  $\cal C$ are isomorphic if and 
only if the angles between the $3$-spaces ${\bf m}_{b_3,c_3,c_2},\  
{\bf m}_{b_3',c_3',c_2'}$ and the $3$-space 
$\hat{{\bf m}}={\bf m}_{0,0,0}$ are the same 
with respect to the Cartan-Killing form.
As representatives of the isomorphism classes of the Bol loops 
$L_{b_3,c_3,c_2}$
may be chosen the Bol loops $\hat{L}_{d}=L_{d,0,0}$  for $d \in [0,1)$.
\newline
The pseudo-euclidean space loop $({\hat L}, \ast )$ has a representation on the 
manifold $\Psi $ of euclidean planes in the $3$-dimensional pseudo-euclidean geometry 
${\cal R}$ with the multiplication is given by the formula $(\ast \ast)$.   
\end{Theo}

\section{Construction of 3-dimensional differentiab\-le 
loops within the pseudo-euclidean space}

Let ${\cal R}$ be the $3$-dimensional pseudo-euclidean affine space. We embed the affine space 
$\mathbb R^3$ into the projective space $P_3(\mathbb R)$ such that $(x,y,z) \mapsto $ 
\newline 
$k (1,x,y,z) $, $k \in \mathbb R 
\backslash \{ 0 \}$. Then the projective plane $E$ at infinity  
consists of the points $\{ k (0,x,y,z) \ |\  x,y,z \in \mathbb R, k \in \mathbb R 
\backslash \{ 0 \} \}$. The cone $x^2+y^2-z^2=0$ 
of the affine space $\mathbb R^3$ intersects the plane $E$ in the conic given by 
$C:\{ k (0,x,y,z)\ |\  x^2+y^2=z^2, k (0,x,y,z) \neq (0,0,0,0) \}$. The conic $C$ divides the 
points of $E \backslash C$ in two regions $R_1, R_2$. The points $(0,x,y,z)$ of $R_1$ 
respectively of $R_2$ are characterized by the property $x^2+y^2< z^2$, respectively
$x^2+y^2>z^2$. 
\newline
Let $\Psi $ be the manifold of planes which intersect the plane $E$ in lines  
contained in the region $R_2$, i.e. in lines which do not meet $C$. 
There is a polarity $\pi $ of $E$ such that the absolute points of $\pi $ are the points of 
the conic $C$. The polarity $\pi $ interchanges the points of $R_1$ with the lines of $R_2$. 
We may assume that the connected component of the motion group  of 
$G=PSL_2(\mathbb R) \ltimes \mathbb R^3$ on $\cal{R}$ is given 
by the Lorentz  transformations. Since any translation of $\mathbb R^3$ 
leaves $E$ pointwise fixed 
the group induced by $G$ on $E$ is isomorphic to $PSL_2(\mathbb R)$ and it is already induced 
by the stabilizer $\Gamma $ of $(0,0,0)$ in $G$.
\newline
Let $\Sigma $ be the submanifold of $\Gamma \cong PSL_2(\mathbb R)$ which is the 
image of the section $\sigma_1$ defining the hyperbolic plane loop (cf. $\sigma _1$ before 
Theorem 6) and let $ \Lambda $ be a $1$-dimensional translation group leaving a line $S$ 
invariant which intersects any plane of $\Psi $ in precisely one point and contains the point 
$(0,0,0)$. We show that the product $\Theta =\Lambda \Sigma$ acts sharply transitively on the 
planes of $\Psi $ and may be taken as the image of a section for a differentiable 
$3$-dimensional loop $L_{\Lambda }$ having $G$ as the group topologically generated by its 
left translations. 
\newline
Let $D_1$, $D_2$ be planes belonging to $\Psi $. We show that there is precisely an element 
$\rho =\omega \delta$ with $\delta \in \Sigma$ and $\omega \in \Lambda$ such that 
$\rho (D_1)=(\omega \delta )(D_1)=D_2$. Let $W_1=D_1 \cap E$ and $W_2=D_2 \cap E$. 
The lines $W_1$ and $W_2$ are contained in the region $R_2$. Hence $\pi (W_1)$ and $\pi (W_2)$ 
are points in $R_1$. In $\Sigma $ there is precisely one element $\delta $ such that 
$\delta (\pi (W_1))=\pi (W_2)$.  But then $\delta (W_1)=W_2$ and $\delta (D_1)$ is parallel to 
$D_2$. Let $s_1= D_1 \cap S$ and $s_2=D_2 \cap S$. Then there exists precisely one translation 
$\omega $ in $\Lambda $ with $\omega (s_1)= s_2$ and $\omega \delta (D_1)=D_2$. 
\newline
Let $P$ be a 
plane of $\Psi $ containing the point $(0,0,0)$ and $H$ be the stabilizer of $P$ in $G$. Then 
for the manifold $\Theta $, which topologically generates $G$, one has $\Theta \cap H=\{ 1 \}$. 
Let $J=P \cap E$ and $L(P)$ be the line joining the points $(0,0,0)$ and $\pi (J)$. This line 
intersects any plane of $\Psi $ in precisely one points. Let $\Lambda _{L(P)}$ be the 
$1$-dimensional translation group leaving the line $L(P)$ invariant. 
\newline
The multiplication for the loop $L_{\Lambda _{L(P)}}$ with the identity $P$ may be defined as 
follows: 
\[ (X,Y) \mapsto \omega_{(\tau _{J,Q}(P),X)} (\tau_{J,Q}(Y)): \Psi \times \Psi \to \Psi ,\]   
where $\tau_{J,Q}$ is the unique element in $\Sigma $ mapping $J$ onto the line 
$Q=X \cap E$ and $\omega _{(\tau _{J,Q}(P),X)}$ is the unique translation in $\Lambda _{L(P)}$ 
mapping the plane $\tau_{J,Q}(P)$ onto $X$. 
\newline
Let $K$ be the $1$-dimensional compact subgroup of $\Gamma $ which leaves the line $J$ 
invariant. Then $K$ fixes $\pi(J)$ and leaves the line $L(P)$ pointwise invariant. Since 
$\Sigma \subset \Gamma $ is the image of a section corresponding to the hyperbolic 
plane loop with the identity $\pi(J)$ we have $g \Sigma g^{-1}= \Sigma$ and  
$g \Lambda_{L(P)} g^{-1}= \Lambda _{L(P)}$ for all $g \in K$. Hence the loop 
$L_{\Lambda _{L(P)}}$ has a $1$-dimensional compact group of automorphisms. 
\newline
Let $\Lambda ^*$ be a $1$-dimensional translation group such that the lines $S$ which are 
invariant under $\Lambda ^*$ intersect any plane of $\Psi $, but $S \cap E \neq \pi(J)$.  
Then the loop $L _{\Lambda ^*}$ with identity $P$ which corresponds to the section 
$\Lambda^* \Sigma $ and to the stabilizer $H$ is not isomorphic to the loop $L_{\Lambda _{L(P)}}$.
This can be seen in the following way: The group $K \subseteq H$ does not leave $S$ invariant 
and hence does not normalize $\Lambda ^*$. It follows that 
$g \Theta g^{-1}=g \Lambda ^* \Sigma g^{-1}=(g \Lambda ^* g^{-1}) g \Sigma g^{-1}=
(g \Lambda ^* g^{-1})\Sigma \neq \Theta$ 
and $K$ does not consists of automorphisms of the loop $L_{\Lambda ^*}$. 
Every loop $L_{g^{-1} \Lambda ^* g}$ with $g \in K$ is isomorphic to $L_{\Lambda ^*}$. 
\newline
No loop $L_{\Lambda }$ is a Bol loop since otherwise for elements $1 \neq \lambda \in \Lambda $ 
and $1 \neq \rho \in \Sigma $ the element $(\lambda \rho) \cdot 1 \cdot \lambda \rho$ should be in 
$\Lambda \Sigma $. But one has $\lambda \rho \lambda \rho=\lambda \rho \lambda \rho^{-1} 
\rho^2 \not\in \Lambda \Sigma$ since $\lambda \rho \lambda \rho^{-1} \not\in \Lambda$. 
The submanifold $\Theta$ is the image of a differentiable section. Hence the multiplication of 
$L_{\Lambda }$ as well as the mappings $(a,b) \mapsto a \backslash b: L_{\Lambda } \times 
L_{\Lambda } \to L_{\Lambda }$ are differentiable.  Let $D_1$ and $D_2$ be elements of 
$\Lambda _L$ and $D_i \cap E=W_i$ for $i=1,2$. Since the mapping $(W_1,W_2) \mapsto 
\tau_{J,Q} \in \Sigma$ with $\tau_{J,Q}(W_1)=W_2$ as well as for parallel planes $D_1$ and 
$D_2$ the map $(D_1, D_2) \to \omega _{D_1, D_2} \in \Lambda$ with $\omega_{D_1, D_2}(D_1)=D_2$ 
are differentiable any loop $L_{\Lambda }$ is differentiable.

\bigskip
\noindent
{\small Mathematisches Institut der Universit\"at Erlangen-N\"urnberg}
\newline
{\small Bismarckstr. 1 $\frac{1}{2}$, D-91054 Erlangen, Germany}
\newline
{\small e-mail: figula@mi.uni-erlangen.de}
\newline
{\small and}
\newline
{\small Institute of Mathematics, University of Debrecen}
\newline
{\small P.O.B. 12, H-4010 Debrecen, Hungary}
\newline
{\small e-mail: figula@math.klte.hu}
\end{document}